\documentclass[leqno]{amsart}
\usepackage{amssymb}
\usepackage{amsmath, amsfonts,enumerate}
\usepackage{color}
\usepackage{amsthm}
\usepackage{hyperref}
\usepackage{mathrsfs}

\newtheorem{theorem}{Theorem}

\newtheorem{proposition}[theorem]{Proposition}
\newtheorem{corollary}[theorem]{Corollary}
\newtheorem{lemma}[theorem]{Lemma}

\newtheorem{definition}[theorem]{Definition}
\theoremstyle{definition}
\numberwithin{theorem}{section}
\numberwithin{equation}{section}

\newtheorem{theoremx}{Theorem}[section]

\def\U{{\mathcal U}}
\def\ri{{\rm{i}}}
\def\sym{{\rm sym}}

\begin{document}

\title[Singular value asymptotics on manifolds]{Singular value asymptotics on compact smooth Riemannian manifolds}



\author{Fedor Sukochev}
\address{School of Mathematics and Statistics, UNSW, Kensington, NSW 2052, Australia}
\email{f.sukochev@unsw.edu.au}

\author{Fulin Yang}
\address{Beijing Institute of Mathematical Sciences and Applications, Beijing 101408, China}
\address{Yau Mathematical Sciences Center, Tsinghua University, Beijing 100084, China}
\email{fulinyang@bimsa.cn}

\author{Dmitriy Zanin}
\address{School of Mathematics and Statistics, UNSW, Kensington, NSW 2052, Australia}
\email{d.zanin@unsw.edu.au}

\date{}

\begin{abstract}
Let $(X,G)$ be a $d$-dimensional compact smooth Riemannian manifold equipped with Laplace-Beltrami operator $\Delta_{G}$, and let $\Pi_{X}$ be the $C^{\ast}$-algebra obtained by locally transferring the $C^{\ast}$-algebra generated by multiplication operators and Riesz transforms on $\mathbb{R}^{d}$. Denote $\sym_{X}$ the principal symbol mapping of $\Pi_{X}$. For any $S\in\Pi_{X}$, we prove that, in the framework of $C^{\ast}$-algebra,
\begin{align*}
\lim_{t\rightarrow\infty}t^{\frac{1}{p}}\mu(t,S(1+\Delta_G)^{-\frac{d}{2p}})
=(2\pi\sqrt[d]{d})^{-\frac{1}{p}}\Big\|\sym_{X}(S)\Big\|_{L_{p}(T^{\ast}X,e^{-q_{G}}d\lambda)}
,
\end{align*}
where $0<p<\infty$, $e^{-q_{G}}$ is the canonical weight on $X$, and $d\lambda$ is the Liouville measure on the cotangent bundle $T^{\ast}X$.
\end{abstract}

\maketitle

\section{Introduction}
In this paper, we establish singular value asymptotics on a $C^{\ast}$-algebra over compact Riemannian manifold.
On $\mathbb{R}^{d}$, recall that $\nabla=(D_{1},\cdots,D_{d})$ and $\Delta=\sum_{k=1}^{d} D_{k}^{2}$, where 
$D_{k}=\frac{1}{i}\partial_{k}$. The $C^{\ast}$-subalgebra $\Pi_{\mathbb{R}^{d}}$ of $B(L_{2}(\mathbb{R}^{d}))$ (the set of bounded linear operators), introduced in \cite{DAO1}, is generated by the algebras
$$\{M_{f}:f\in\mathbb{C}+C_{0}(\mathbb{R}^{d})\},\quad
\{g(\frac{\nabla}{\sqrt{\Delta}}):g\in C(\mathbb{S}^{d-1})\},$$
where $M_{f}$ is defined by $M_{f}\phi=f\phi$ for $\phi\in L_{1,\text{loc}}(\mathbb{R}^{d})$. 
As shown in \cite{DAO1,DAO2}, there exists a $\ast$-homomorphism
$${\rm sym}:\Pi_{\mathbb{R}^{d}}\rightarrow C(\mathbb{S}^{d-1},\mathbb{C}+C_{0}(\mathbb{R}^{d}))$$
such that, for all $f\in\mathbb{C}+C_{0}(\mathbb{R}^{d})$ and $g\in C(\mathbb{S}^{d-1})$,
$${\rm sym}(M_{f})=f\otimes1\quad\mbox{and}\quad{\rm sym}(g(\frac{\nabla}{\sqrt{\Delta}}))=1\otimes g.$$
This $*$-homomorphism, referred to as a principal symbol mapping, generalizes the classical notion of the principal symbol for pseudodifferential operators on $\mathbb{R}^{d}$.

An operator $S\in B(L_{2}(\mathbb{R}^{d}))$ is compactly supported if there exists $\phi\in C_{c}^{\infty}(\mathbb{R}^{d})$ such that $S=M_{\phi}S M_{\phi}$. For a compact operator $T$ on Hilbert space $H$, let $t\mapsto\mu(t,T)$ be the singular value function (see Subsection \ref{Cwikel estimates} for detail).
The following singular value asymptotics on $\mathbb{R}^{d}$ was proved in\cite[Theorem 6.1]{FSZ2024endpoint} 
(earlier in \cite{FSZ2023asymptotics} with the case $p=d$).
\begin{theoremx}\label{asymptotic on Rd}
For $0 < p < \infty$, if $S \in \Pi_{\mathbb{R}^{d}}$ is compactly supported, then
$$\lim_{t\rightarrow\infty} t^{\frac{1}{p}}\mu\Big(t,S(1+\Delta)^{-\frac{d}{2p}}\Big)
=c_d^{\frac1p}\|\sym(S)\|_{L_p(\mathbb{R}^{d}\times\mathbb{S}^{d-1})}$$
where $c_{d}=d^{-\frac{1}{d}}(2\pi)^{-1}$.
\end{theoremx}
The Theorem \ref{asymptotic on Rd} is closely connected to the Connes' integral established in \cite[Theorem 5.3.3]{LMSZvol2}(also \cite[Theorem 1.5]{DAO1}), and plays a key role in studying the singular value asymptotics of quantum derivatives \cite{FSZ2023asymptotics, SXZ2023, Ti2025} and spectral estimates \cite{FSZ2024endpoint,FLSZ2024}.

Let $(X, G)$ be a compact smooth Riemannian manifold with Riemannian metric $G$ and an atlas $\{(\mathcal{U}_{i},h_{i})\}_{i\in\mathbb{I}}$, and let ${\rm vol}_{G}$ be the volume form on $X$ induced by the metric $G$. On the compact manifold $X$, the $C^{\ast}$-algebra $\Pi_{X}$ introduced in \cite{DAO3} corresponds to $\Pi_{\mathbb{R}^{d}}$ in local coordinates.
\begin{definition}
Let $(X,G)$ be a $d$-dimensional compact smooth Riemannian manifold with Riemannian metric $G$. For $S \in B(L_{2}(X,{\rm vol}_{G}))$, we say $S\in\Pi_{X}$ if
\begin{enumerate}[\rm(1)]
\item For each $i\in\mathbb{I}$ and each $\phi\in C_{c}(\mathcal{U}_{i})$, the operator 
$M_{\phi}SM_{\phi}$, when transferred to an operator on $L_{2}(\mathbb{R}^{d})$, belongs to $\Pi_{\mathbb{R}^{d}}$;
\item For each $\psi\in C(X)$, the commutator $[M_{\psi},S]$ is compact.
\end{enumerate}
\end{definition}
The $C^{\ast}$-algebra $\Pi_{\mathbb{R}^{d}}$ is strictly broader than the class of classical compactly supported pseudo-differential operators of order 0 on $\mathbb{R}^{d}$ \cite{KN1965}.  By the definition, the algebra $\Pi_{X}$ contains the class of classical compactly supported pseudo-differential operators of order 0 in each chart. 

Let $S^{\ast}X$ denote the cosphere bundle of $X$. By \cite[Theorem 1.4]{DAO3}, the algebra $\Pi_{X}$ is a $C^{\ast}$-algebra, and there exists a surjective $\ast$-homomorphism
$${\rm sym}_{X}:\Pi_{X}\rightarrow C(S^{*}X)$$ 
such that
$${\rm ker}({\rm sym}_{X})=\mathcal{K}(L_{2}(X,{\rm vol}_{G})).$$
Here, $\mathcal{K}(L_{2}(X,{\rm vol}_{G}))$ is the set of compact operators on Hilbert space $L_{2}(X,{\rm vol}_{G})$.
We prove an analogue of Theorem \ref{asymptotic on Rd} on compact Riemannian manifolds.
\begin{theorem}\label{asymptotic on manifold}
Let $(X,G)$ be a compact smooth Riemannian manifold with Riemannian metric $G$. For $0<p<\infty$, if $S\in\Pi_{X}$, then
$$\lim_{t\rightarrow\infty}t^{\frac{1}{p}}\mu(t,S(1+\Delta_G)^{-\frac{d}{2p}})
=(2\pi\sqrt[d]{d})^{-\frac{1}{p}}\Big\|\sym_{X}(S)\Big\|_{L_{p}(T^{\ast}X,e^{-q_{G}}d\lambda)}$$
where $\Delta_{G}$ and $e^{-q_{G}}$ are the Laplace–Beltrami operator and canonical weight on $X$ separately, and $d\lambda$ is the Liouville measure on the cotangent bundle $T^{\ast}X$.
\end{theorem}
For the precise definitions of the Laplace–Beltrami operator $\Delta_{G}$, the weight $e^{-q_{G}}$, and the Liouville measure $d\lambda$, we refer to Subsection \ref{def of compact manifold}.

The Theorem \ref{asymptotic on manifold} covers the results on closed Riemannian manifolds \cite{BS1977,Gru1984} and the classical tori \cite{SXZ2023} (this result is different with that in \cite{BS1977,Gru1984}), and has a direct link with Connes' integral in \cite{DAO2,DAO3} on compact Riemannian manifolds. Let us give the details below.

\subsection{Motivation on Weyl's law}
For a bounded domain $\Omega\subset\mathbb{R}^{d}$ with Dirichlet boundary conditions, H. Weyl \cite{Weyl1912} found that the eigenvalue counting function $N(\lambda)$ for the Laplacian $\Delta$ below $\lambda$ satisfies the asymptotic formula:
$$N(\lambda)=\frac{\omega_{d}{\rm vol}(\Omega)}{(2\pi)^{d}}\lambda^{\frac{d}{2}}+o(\lambda^{\frac{d}{2}}),
\quad\mbox{as}\;\lambda\rightarrow\infty,$$
where $\omega_{d}=\frac{\pi^{\frac{d}{2}}}{\Gamma(\frac{d}{2}+1)}$ is the volume of the unit ball in $\mathbb{R}^{d}$, and ${\rm vol}(\Omega)$ is the volume of $\Omega$. This result, now known as Weyl's law, was later extended by Birman--Solomyak \cite{BS1977} to negative-order pseudo-differential operators in local settings of $\mathbb{R}^d$, and generalized to compact Riemannian manifolds by Chazarain \cite{Ch1974} (see also \cite[page 155]{Ch1984}). 
In recent years, Weyl's law has been actively investigated from both geometric \cite{Po2008,GG2019,Fl2024,CPR2024} and $C^{\ast}$-algebraic \cite{MSZ2022,FSZ2023asymptotics,SXZ2023,FSZ2024endpoint} perspectives, with ongoing advances in results and associated techniques.

A powerful modern reformulation of Weyl's law is expressed in terms of singular value asymptotics. Let $X$ be a $d$-dimensional closed Riemannian manifold, and let $T$ be a pseudo-differential operator of order $-p$ on $X$ with $p>0$. In \cite[Lemma 4.5]{Gru1984}, Grubb states (she referred to \cite{BS1977} without proof) that
\begin{align}\label{classical asymptotics}
\lim_{k\rightarrow\infty}k^{\frac{p}{d}}\mu(k,T)=
\frac{1}{d(2\pi)^{d}}\|\sigma_{T}\|_{L_{\frac{d}{p}}(S^{*}X)}
\end{align}
where $\sigma_{T}$ is the principal symbol of $T$ and $S^{\ast}X$ is the cosphere bundle of $X$. 
The formula \eqref{classical asymptotics} is related to Connes's original noncommutative integration \cite{Co1994,Ro2022,Po2023}. It is generalised to Theorem \ref{asymptotic on Rd} on Euclidean space \cite{FSZ2023asymptotics,FSZ2024endpoint} and a counterpart of Theorem \ref{asymptotic on Rd} (with the case $p=d$) on quantum tori \cite{SXZ2023} including the classical tori. 

\subsection{Connection with noncommutative integration} 
Let $(X,G)$ be a compact smooth Riemannian manifold of dimension $d$, equipped with Laplace-Beltrami operator $\Delta_{G}$. In noncommutative geometry \cite{Co1994}, A. Connes observed that Weyl's law could be reinterpreted as follows: the integral of $f\in C^{\infty}(X)$ equals the Dixmier trace of operator $M_{f}(1+\Delta_{G})^{-\frac{d}{2}}$.
In \cite{Co1994}, the functional 
$${\rm Int}(T):A\mapsto{\rm Tr}_{\omega}\left(T(1+\Delta_{G})^{-\frac{d}{2}}\right),\quad T\in B(L_{2}(X,{\rm vol}_{G})),$$ 
where ${\rm Tr}_{\omega}$ denote the Dixmier trace with respect to an extended limit $\omega$ on $\ell_{\infty}$, is proposed as a notion of noncommutative integration. Whenever $T(1+\Delta_{G})^{-\frac{d}{2}}$ belongs to $\mathcal{L}_{1,\infty}(L_{2}(X,{\rm vol}_{G}))$, this definition is independent of the choice of the extended limit $\omega$ and yields the same value  (see \cite{LMSZvol2} for more details). The connection between asymptotic formulas for singular values and noncommutative integral formulas was first studied by Rozenblum \cite{Ro2022}, and later explored in \cite{SZ-FAA}.

According to \cite[Theorem 1.5]{DAO3} (see also \cite{DAO2,LMSZvol2} for the case of classical tori), for any $S\in\Pi_{X}$, one has $T(1+\Delta_{G})^{-\frac{d}{2}}\in\mathcal{L}_{1,\infty}(L_{2}(X,{\rm vol}_{G}))$ and
$${\rm Int}(S)=c_{d}\int_{T^{\ast}X}{\rm sym}_{X}(S)e^{-q_{G}}d\lambda.$$
If $S\in\Pi_{X}$ is positive, Theorem \ref{asymptotic on manifold} states that 
$${\rm Int}(S)=\lim_{t\rightarrow\infty}t\mu(t,S(1+\Delta_G)^{-\frac{d}{2}}).$$
Hence, Theorem \ref{asymptotic on manifold} shows a direct and explicit link between noncommutative integral formula in \cite{DAO3} and singular value asymptotics for positive elements in $\Pi_{X}$.

\subsection{Advantages and obstacles} 
In contrast to \eqref{classical asymptotics}, the advantage of Theorem \ref{asymptotic on manifold} is twofold: 
(a) it encompasses a strictly larger class of operators; 
(b) it is situated within the category of $C^{\ast}$-algebras (i.e. non-commutative topological spaces), not in a category of classical pseudodifferential operators which does not have a natural counterpart in non-commutative geometry. In addition, Theorem \ref{asymptotic on manifold} covers a counterpart of Theorem \ref{asymptotic on Rd} on classical tori \cite{SXZ2023} and provides a link with \cite{DAO3}  Connes' noncommutative integral for positive operators in $\Pi_{X}$.

Several significant obstacles arise to treat Theorem \ref{asymptotic on manifold}. First, since the Laplace--Beltrami operator locally equals to an elliptic differential operator, we need to establish spectral estimates for commutators involving elliptic differential operators on $\mathbb{R}^{d}$. Second, we must adapt Theorem \ref{asymptotic on Rd} to a curved space setting with respect to an elliptic differential operator on $\mathbb{R}^{d}$. Third, transferring singular value asymptotics from local to global settings requires avoiding the Complex Powers Theorem for pseudo-differential operators on manifolds, as no satisfactory version is available.

Our proof diverges from the methods of Birman and Solomyak \cite{BS1977,Gru1984} by employing a novel combination of functional calculus and double operator integrals. Two key innovations are the differentiability of operators inspiring by Connes--Moscovici (Subsection \ref{Application of abstract estimate}) and the new concept of ``good atlas" on compact manifolds (Subsection \ref{Good atlas}). We demonstrate the density of differentiable elements in $\Pi_{X}$ (Lemma \ref{main example density lemma}) and the existence of ``good atlas" for any compact Riemannian manifold (Proposition \ref{good atlas prop}), which serve as foundations to show Theorem \ref{asymptotic on manifold}. Collectively, these strategies go beyond the classical pseudo-differential calculus.

\vspace{0.5em}
\subsection{Outline of the Paper}

Section \ref{sec-2} provides necessary notions and facts.

In Section \ref{sec-3}, we establish spectral estimates (Theorems \ref{euclidean compactness of commutators} and \ref{1 over 2 case}) for commutators of elliptic differential operators on $\mathbb{R}^{d}$ arising from the Laplace–Beltrami operator on manifolds.

In Section \ref{sec-4}, we generalize Theorem \ref{asymptotic on Rd} to Theorem \ref{for curved plane thm} by replacing the Laplacian with an elliptic differential operator on $\mathbb{R}^{d}$ for the index $p = \frac{d}{2}$.

In Section \ref{sec-5}, we give spectral characterization (Theorem \ref{main example thm}) for powers of two nearly commutating operators.

In Section \ref{sec-6}, we show Theorem \ref{asymptotic on manifold} in a local chart by using the spectral estimates from Section \ref{sec-3} and the singular value asymptotics from Section \ref{sec-4}.

Finally, in Section \ref{sec-7}, we combine the spectral estimates from Sections \ref{sec-3} and \ref{sec-5} with the local version of Theorem \ref{asymptotic on manifold} from Section \ref{sec-6} to end the proof.

Throughout, we use $c_A$ to denote some  constant depending only on the symbol $A$,  which may change from line to line without specification.

\section{Preliminaries}\label{sec-2}
In this section, we recall necessary notions for manifold and double operator integral, and provide some facts related to Cwikel estimate and asymptotics.
\subsection{Compact manifold}\label{def of compact manifold}
The following materials are standard, see e.g. \cite{Ch2006,Mu2018}. Let $(X,G)$ be a $d$-dimensional compact smooth Riemannian manifold with Riemannian metric $G$ and an atlas $\{(\U_{i},h_{i})\}_{i\in\mathbb{I}}$. 
Let $$\Omega_{i}=h_{i}(\U_{i})\subset\mathbb{R}^{d},\quad
\Omega_{i,j}=h_{i}(\U_{i}\cap\U_{j})\subset\mathbb{R}^{d}.$$
Denote $\Phi_{i,j}:\Omega_{i,j}\rightarrow\Omega_{j,i}$ the diffeomorphism given by the formula
$\Phi_{i,j}=h_{j}\circ h_{i}^{-1}.$
Let $T^{\ast}X$ be the cotangent bundle of $X$ and let $\pi:T^{\ast}X\rightarrow X$ be the canonical projection. There exists an atlas $\{(\pi^{-1}(\U_{i}),H_{i})\}_{i\in\mathbb{I}}$ on $T^{\ast}X$ such that
\begin{enumerate}[\rm(a)]
\item for $i\in\mathbb{I}$, function $H_{i}:\pi^{-1}(\U_{i})\rightarrow\Omega_{i}\times\mathbb{R}^{d}$ is a homeomorphism;
\item for $i,j\in\mathbb{I}$ such that $\Omega_{i,j}\neq\emptyset$, we have 
$$(H_j\circ H_i^{-1})(\xi,s)
=(\Phi_{i,j}(\xi),(J_{\Phi_{i,j}}^{\ast}(\xi))^{-1}s),\quad\xi\in\Omega_{i,j},\quad s\in\mathbb{R}^{d},$$
where $J_{\Phi_{i,j}}$ is the Jacobian matrix of mapping $\Phi_{i,j}$.
\end{enumerate}
Since $T^{\ast}X$ is a smooth $d$-dimensional manifold, it follows that $T^{\ast}X$ has a canonical symplectic structure \cite{GS1990}. The corresponding Liouville measure $d\lambda$ on $T^{\ast}X$ is given by 
$$\int_{T^{\ast}X}(f\circ H_{i})d\lambda=\int_{\mathbb{R}^{d}\times\mathbb{R}^{d}}f(\xi,s)dm(\xi, s),\quad f\in C_{c}(\Omega_{i}\times\mathbb{R}^{d}),\;i\in\mathbb{I},$$
where $dm$ is the product-Lebesgue measure on $\mathbb{R}^{d}\times\mathbb{R}^{d}$.
Here $f\circ H_{i}$ denotes a function on $T^{\ast}X$ which equals $f\circ H_{i}$ on $\pi^{-1}(\U_{i})$ and which vanishes outside $\pi^{-1}(\U_{i})$. By the way, there is no canonical way to equip the cosphere bundle $S^{\ast}X$ of a smooth manifold $X$ with a measure.

Let $\Omega\subset\mathbb{R}^{d}$ be connected and open. Let $g:\Omega\rightarrow{\rm GL}^{+}(d,\mathbb{R})$ be a smooth mapping where ${\rm GL}^{+}(d,\mathbb{R})$ stands for the set of all positive elements in ${\rm GL}(d,\mathbb{R})$ (the space of all $d\times d$ invertible real matrices). The Laplace-Beltrami operator $\Delta_{g}:C_{c}^{\infty}(\Omega)\rightarrow C_{c}^{\infty}(\Omega)$ associated to mapping $g$ is defined by 
$$\Delta_{g}=M_{{\rm det}(g)^{-\frac{1}{2}}}\sum_{k,l=1}^{d}D_{k}M_{{\rm det}(g)^{\frac{1}{2}}\cdot(g^{-1})_{kl}}D_{l}$$
where $D_{k}=\frac1i\partial_{k}$ and ${\rm det}(g)$ is the determinant of $g$. 

For $i\in\mathbb{I}$, the components of the metric $G$ in the chart $(\U_{i},h_{i})$ give rise to a smooth mapping $G_{i}:\U_{i}\rightarrow{\rm GL}^{+}(d,\mathbb{R})$. 
Set $g_{i}=G_{i}\circ h_{i}^{-1}:\Omega_{i}\rightarrow{\rm GL}^{+}(d,\mathbb{R})$. Define a partial isometry $W_{h_i}:L_2(X,{\rm vol}_G)\to L_2(\mathbb{R}^d,{\rm vol}_{g_i})$ by setting 
$$W_{h_i}\xi=\xi\circ h_i^{-1}.$$	
Clearly,
$$W_{h_i}^{\ast}W_{h_i}=M_{\chi_{\U_i}},\quad W_{h_i}W_{h_i}^{\ast}=M_{\chi_{h_i(\U_i)}}.$$

Since $X$ is compact, there is a good partition of unity $\{\phi_{n}\}_{n=1}^{N}\subset C_{c}^{\infty}(X)$ which is a finite partition of unity on $X$ and each $\phi_{n}$ is compactly supported in some chart.
The Laplace-Beltrami operator $\Delta_{G}:C^{\infty}(X)\rightarrow C^{\infty}(X)$ is defined by the formula $$\Delta_{G}f=\sum_{n=1}^{N}\Delta_{g_{i_{n}}}\Big((f\phi_{n})\circ h_{i_{n}}^{-1}\Big)\circ h_{i_{n}},\quad f\in C^{\infty}(X).$$
Here $i_{n}$ is chosen such that $\phi_{n}$ is compactly supported in $\U_{i_{n}}$. 
This definition does not rely on the choice of the good partition of unity $\{\phi_{n}\}_{n=1}^{N}$.

For our purpose, we need to introduce the canonical weight $e^{-q_{G}}$ on $T^{\ast}X$. Set $$q_{i}(\xi,s)=\langle g_{i}(\xi)^{-1}s,s\rangle,\quad\xi\in\Omega_{i},\;s\in\mathbb{R}^{d}.$$
By \cite[Fact 2.10]{DAO3}, there exists a function $q_{G}$ defined on $T^{\ast}X$ such that
$$q_{i}=q_{G}\circ H_{i}^{-1}\quad\mbox{on}\quad\Omega_{i}\times\mathbb{R}^{d}.$$ 
The function $q_{G}$ is the square of the length function on $T^{\ast}X$ defined by the induced Riemannian metric on the cotangent bundle $T^{\ast}X$. More exactly, in local coordinates, the function $e^{-q_{G}}$ satisfies $e^{-q_{G}}\circ H_{i}^{-1}(\xi,s)=e^{-q_{i}(\xi,s)}$ on $\Omega_{i}\times\mathbb{R}^{d}$. We call  $e^{-q_{G}}$ the canonical weight of $(X,G)$ on $T^{\ast}X$.
Since $X$ is compact, it follows $e^{-q_{G}}\in L_{1}(T^{\ast}X,d\lambda)$.

In what follows, we always suppose that the manifold in considering is smooth. 

\subsection{Pseudo-differential operator}
For $a\in C^{\infty}(\mathbb{R}^{d}\times\mathbb{R}^{d})$, define
$${\rm Op}(a)f(x)=\int_{\mathbb{R}^{d}}e^{2\pi{\rm i}\langle x,\xi\rangle}a(x,\xi)\hat{f}(\xi)d\xi,\quad f\in L_{2}(\mathbb{R}^{d}),$$
where $\hat{f}$ is the Fourier transform of $f$. By the Calder\'{o}n-Vaillancourt theorem, ${\rm Op}(a)$ is a bounded on $L_{2}(\mathbb{R}^{d})$.
For $m\in\mathbb{Z}$, let $\Psi^{m}(\mathbb{R}^{d})$ be the class of operators ${\rm Op}(a)$ with its symbol $a$ satisfying 
$$|D_{x}^{\alpha}D_{\xi}^{\beta}a(x,\xi)|\leq c_{\alpha,\beta}(1+|\xi|)^{m-|\beta|_{1}},
\quad \alpha,\beta\in\mathbb{Z}_{+}^{d}.$$
Here $|\beta|_{1}=\sum_{k=1}^{d}|\beta_{k}|$ and $|\xi|=(\sum_{k=1}^{d}|\xi_{k}|^{2})^{\frac{1}{2}}$. The elements in $\Psi^{m}(\mathbb{R}^{d})$ are called pseudo-differential operators of order $m$. The key properties are that 
\begin{align}\label{pseudo-differential op}
\Psi^{m_{1}}(\mathbb{R}^{d})\cdot \Psi^{m_{2}}(\mathbb{R}^{d})\subset\Psi^{m_{1}+m_{2}}(\mathbb{R}^{d}),\quad [\Psi^{m_{1}}(\mathbb{R}^{d}),\Psi^{m_{2}}(\mathbb{R}^{d})]\subset\Psi^{m_{1}+m_{2}-1}(\mathbb{R}^{d})
\end{align}
for $m_{1},m_{2}\in\mathbb{Z}$.

Let $(X,G)$ be a compact Riemannian manifold. A linear operator $A:C^{\infty}(X)\rightarrow C^{\infty}(X)$ is called a pseudo-differential operator on $X$ of order $m$ if, for any chart $(\mathcal{U},h)$ and any $\psi,\phi\in C_{c}^{\infty}(\mathcal{U})$, the operator $M_{\psi}AM_{\phi}:C^{\infty}(\mathcal{U})\rightarrow C^{\infty}(\mathcal{U})$ is a pseudo-differential operator on $\mathbb{R}^{d}$ of order $m$. Denote $\Psi^{m}(X)$ the class of all pseudo-differential operators on $X$ of order $m$. Clearly, $\Delta_{G}\in \Psi^{2}(X)$. Similarly,
$$\Psi^{m_{1}}(X)\cdot \Psi^{m_{2}}(X)\subset\Psi^{m_{1}+m_{2}}(X),\quad 
[\Psi^{m_{1}}(X),\Psi^{m_{2}}(X)]\subset\Psi^{m_{1}+m_{2}-1}(X)$$
for $m_{1},m_{2}\in\mathbb{Z}$. 

More information on pseudo-differential operators, we refer \cite{DAO3, Stein-Murphy1993} for $\mathbb{R}^{d}$ and \cite{Shubin-book2001,RT2009-book} for compact manifolds.

\subsection{Sobolev space}
The Sobolev space $W^{m,2}(\mathbb{R}^{d})$, $m\in\mathbb{N}$, consists of all distributions $f\in L_{2}(\mathbb{R}^{d})$ such that distribution $(1+\Delta)^{\frac{m}{2}}f$ belongs to $L_{2}(\mathbb{R}^{d})$.

A linear differential operator $$P=\sum_{k,l=1}^dM_{a_{kl}}D_kD_l+\sum_{k=1}^dM_{b_k}D_k+M_f,\quad a_{kl},b_{k},f\in L_{1,loc}(\mathbb{R}^{d})$$
is said to be elliptic if the matrix $(a_{kl}(x))_{k,l=1}^{d}$ is positive definite for all $x\in\mathbb{R}^{d}$. The Sobolev spaces $W^{2,2}(\mathbb{R}^{d})$ is the domain of the positive operator $\Delta$. Moreover, we crucially use the following fundamental result related to elliptic differential operator which is organised from \cite{Shubin-book2001}. 
\begin{theorem}\label{standard elliptic theorem} 
Let $P:W^{2,2}(\mathbb{R}^d)\to L_2(\mathbb{R}^d)$ be a positive, self-adjoint, elliptic differential operator whose coefficients are constant outside some ball. For every $\alpha>0,$ we have $(1+P)^{-\frac{\alpha}{2}}:L_2(\mathbb{R}^d)\to W^{\alpha,2}(\mathbb{R}^d)$ is bounded.
\end{theorem}

Let $(X,G)$ be a compact Riemannian manifold and let $\Delta_{G}$ be the related Laplace-Beltrami operator. For any $m\in\mathbb{N}$, a distribution $f\in L_{2}(X,{\rm vol}_{G})$ belongs to the Sobolev space 
$W^{m,2}(X,{\rm vol}_{G})$ if and only if $$(1+\Delta_{G})^{\frac{m}{2}}f\in L_{2}(X,{\rm vol}_{G})$$
where the operator is defined spectrally. This global, operator-theoretic characterization is equivalent to the local coordinate-based definition via a good partition of unity, a consequence of the manifold's compactness. If $P:C^{\infty}(X)\rightarrow C^{\infty}(X)$ is a pseudo-differential operator of order $m$, then $P$ has a continuous extension $P:W^{k,2}(X,{\rm vol}_{G})\rightarrow W^{k-m,2}(X,{\rm vol}_{G})$.

\subsection{Cwikel estimates}\label{Cwikel estimates}
Given a Hilbert space $H$, denote $B(H)$ the algebra of all bounded linear operators on $H$. The ideal of all compact linear operators is denoted by $\mathcal{K}(H)$. Given $T\in\mathcal{K}(H)$, the singular value sequence $\mu(T)=\{\mu(n,T)\}_{n=0}^{\infty}$ is defined in term of the singular value function
$$\mu(t,T)=\inf\{\|T-R\|_{B(H)}:{\rm Rank}(R)\leq t\},\quad t\geq0.$$
When necessary, we write $\mu_{B(H)}(t, T)$ to emphasize the ambient algebra. 

For $0<p<\infty$, the Schatten ideal $\mathcal{L}_{p}(H)$ consists of all compact operators $T$ such that $\mu(T)\in\ell_{p}$. In this paper, we concern on the weak Schatten ideal $\mathcal{L}_{p,\infty}(H)$, which consists of all compact operators $T$ with $\mu(T)\in\ell_{p,\infty}$ and (quasi)-norm $$\|T\|_{\mathcal{L}_{p,\infty}(H)}:=\sup_{n\geq0}(n+1)^{\frac{1}{p}}\mu(n,T).$$
For convenience, $B(H)$ is also denoted by $\mathcal{L}_{\infty}(H)$. One has the following H\"{o}lder-type inequality
$$\|TS\|_{\mathcal{L}_{r,\infty}(H)}\leq c_{p,q,r}
\|T\|_{\mathcal{L}_{p,\infty}(H)}\|T\|_{\mathcal{L}_{q,\infty}(H)},\quad\frac{1}{r}=\frac{1}{p}+\frac{1}{q}
,\quad p,q,r\geq1.$$

For the purpose, denote by $(\mathcal{L}_{p,\infty})_{0}(H)$ the separable part of $\mathcal{L}_{p,\infty}(H)$, defined as the subset of $T\in\mathcal{L}_{p,\infty}(H)$ such that 
$$\lim_{n\rightarrow\infty}(n+1)^{\frac{1}{p}}\mu(n,T)=0.$$
A fact we will use repeatedly is that $\mathcal{L}_{p,\infty}(H)\subset(\mathcal{L}_{p+\epsilon,\infty})_{0}(H)$ for any $\epsilon>0$. For further details on trace ideals and the singular value function, see e.g. \cite{IM-1969,Simon2005,LMSZvol2}.

Let us give some relevant Cwikel estimates here on $\mathbb{R}^{d}$ and compact manifold.
\begin{lemma}\label{symmetric Cwikel estimate}
Let $\alpha,\gamma\geq0$ be such that $\alpha+\gamma>0$. If $f\in C_{c}^{\infty}(\mathbb{R}^{d})$, then
$$(1+\Delta)^{-\frac{\alpha}{2}}M_{f}(1+\Delta)^{-\frac{\gamma}{2}}
\in\mathcal{L}_{\frac{d}{\alpha+\gamma},\infty}(L_{2}(\mathbb{R}^{d})).$$
\end{lemma}
\begin{proof}
Write $f=f_{1}f_{2}$ with $f_{1},f_{2}\in C_{c}^{\infty}(\mathbb{R}^{d})$. By H\"{o}lder's inequality, it suffices to show 
$$(1+\Delta)^{-\frac{\alpha}{2}}M_{f_{1}}
\in\mathcal{L}_{\frac{d}{\alpha},\infty}(L_{2}(\mathbb{R}^{d})),\quad
M_{f_{2}}(1+\Delta)^{-\frac{\gamma}{2}}
\in\mathcal{L}_{\frac{d}{\gamma},\infty}(L_{2}(\mathbb{R}^{d})).$$
Here $\mathcal{L}_{\frac{d}{\alpha},\infty}(L_{2}(\mathbb{R}^{d}))$ ($\mathcal{L}_{\frac{d}{\gamma},\infty}(L_{2}(\mathbb{R}^{d}))$) is understood as $\mathcal{L}_{\infty}(L_{2}(\mathbb{R}^{d}))$ if $\alpha=0$ ($\gamma=0$). The inclusions follow from \cite[Theorem 2.4]{FSZ2024endpoint}. 
\end{proof}

\begin{lemma}\label{special cwikel manifold} 
Let $(X,G)$ be a $d$-dimensional compact Riemannian manifold. For well-defined function $\psi$ on $X$, we have
$$\|M_{\psi}(1+\Delta_{G})^{-\frac{d}{2p}}\|_{\mathcal{L}_{p,\infty}(L_{2}(X,{\rm vol}_{G}))}\leq c_{p,X,G}\|\psi\|_{L_p(X,{\rm vol}_{G})},\quad p>2,$$
$$\|M_{\psi}(1+\Delta_{G})^{-\frac{d}{2p}}\|_{\mathcal{L}_{p,\infty}(L_{2}(X,{\rm vol}_{G}))}\leq c_{p,X,G}\|\psi\|_{L_2(X,{\rm vol}_{G})},\quad 0<p<2,$$
$$\|M_{|\psi|^{\frac{1}{2}}}(1+\Delta_{G})^{-\frac{d}{2p}}\|_{\mathcal{L}_{p,\infty}(L_{2}(X,{\rm vol}_{G}))}\leq c_{p,X,G}\|\psi\|^{\frac{1}{2}}_{L_{\Phi}(X,{\rm vol}_{G})},\quad p=2,$$
where $L_{\Phi}(X,{\rm vol}_{G})$ is the Orlicz space associated with the Young function $\Phi(t)= t\log(e+t)$, $t>0$.
\end{lemma}
\begin{proof} The hardest case $p=2$ is established in \cite[Theorem 1.3]{SZ-FAA}. The proof in cases $p>2$ and $p<2$ is virtually identical. One simply needs to replace special Cwikel estimates established in \cite{SZ-Sb} with (a) Cwikel estimates for $p>2$ (b) Birman-Solomyak version of Cwikel estimates for $p<2$ --- see \cite{LSZ-PLMS}.
\end{proof}

\subsection{Birman-Solomyak lemmas}
Fix $H$ a separable Hilbert space. The operators $\{A_{k}\}_{k\geq0}\subset B(H)$ are called pairwise orthogonal if $A_kA_l=A_k^{\ast}A_l=A_kA_l^{\ast}=0$ whenever $k\neq l.$ Let us recall the following lemmas.
\begin{lemma}\cite[Lemma 3.3]{FSZ2023asymptotics}\label{bs direct sum lemma}
Let $\{A_k\}_{k=1}^K\subset\mathcal{L}_{p,\infty}(H)$ be a sequence of pairwise orthogonal operators such that, for all $1\leq k\leq N$, the limits
$$\lim_{t\rightarrow\infty}t^{\frac{1}{p}}\mu(t,A_{k})=c_{k}\quad\mbox{exist}.$$
Then $$\lim_{t\rightarrow\infty}t^{\frac{1}{p}}\mu(t,\sum_{k=1}^KA_{k})=(\sum_{k=1}^Kc_{k}^{p})^{\frac{1}{p}}.$$
\end{lemma}

An important feature for asymptotics of $(\mathcal{L}_{p,\infty})_{0}$ is the following.
\begin{lemma}\label{bs separable lemma} 
Let $0<p<\infty.$ Let $A\in\mathcal{L}_{p,\infty}(H)$ be such that 
$$\lim_{t\to\infty}t^{\frac1p}\mu(t,A)=c.$$
For every $B\in(\mathcal{L}_{p,\infty})_0(H),$ we have
$$\lim_{t\to\infty}t^{\frac1p}\mu(t,A+B)=c.$$
\end{lemma}
The elementary fact is ideed due to K. Fan, see e.g. \cite[Chapter II, Theorem 4.1]{IM-1969}. 
Its a generalisation due to Birman and Solomyak is the following one.

\begin{lemma}\cite{RSS1989,FSZ2023asymptotics}\label{bs limit lemma}
Let $0<p<\infty$. Let $\{A_{n}\}_{n\geq1}\subset\mathcal{L}_{p,\infty}(H)$ be such that
\begin{enumerate}[\rm(a)]
\item $A_{n}\rightarrow A$ in $\mathcal{L}_{p,\infty}(H)$ as $n\rightarrow\infty$;
\item for each $n\geq1$, the limit $\displaystyle\lim_{t\rightarrow\infty}t^{\frac{1}{p}}\mu(t,A_{n})=c_{n}$ exists.
\end{enumerate}
Then the following limits exist and are equal,
$$\lim_{t\rightarrow\infty}t^{\frac{1}{p}}\mu(t,A)=\lim_{n\rightarrow\infty}c_{n}.$$
\end{lemma}
These three lemmas will play important roles in the proof of Theorem \ref{asymptotic on manifold}.

\section{Compactness of commutator on Euclidean space}\label{sec-3}
In this section, we prove Theorems \ref{euclidean compactness of commutators} and \ref{1 over 2 case}, the spectral estimates related to elliptic differential operator on $\mathbb{R}^{d}$ arising from the Laplace–Beltrami operator on manifold.

\begin{theorem}\label{euclidean compactness of commutators}
Let $P:W^{2,2}(\mathbb{R}^d)\to L_2(\mathbb{R}^d)$ be a positive, self-adjoint, elliptic differential operator whose coefficients are constant outside some ball. Let $\phi\in C^{\infty}_c(\mathbb{R}^d)$. If $\alpha,\beta,\gamma\in\mathbb{R}$ satisfy $\alpha+\beta+\gamma+1>0$, then
$$(1+P)^{-\frac{\alpha}{2}}[M_{\phi},(1+P)^{-\frac{\beta}{2}}](1+P)^{-\frac{\gamma}{2}}
\in\mathcal{L}_{\frac{d}{\alpha+\gamma+\beta+1},\infty}(L_2(\mathbb{R}^d)).$$
\end{theorem}

\begin{theorem}\label{1 over 2 case}
Let $P:W^{2,2}(\mathbb{R}^d)\to L_2(\mathbb{R}^d)$ be a positive, self-adjoint, elliptic differential operator whose coefficients are constant outside some ball. Let $0<s<2$ and $0\leq\alpha,\gamma<\infty$ with $\alpha+\gamma+s+1>0$. If $\phi\in C_{c}^{\infty}(\mathbb{R}^{d})$, then 
\begin{align*}
(1+P)^{-\frac{\alpha}{2}}[M_{\phi}R_{j}M_{\phi},(1+P)^{-\frac{s}{2}}](1+P)^{-\frac{\gamma}{2}}
\in\mathcal{L}_{\frac{d}{\alpha+\gamma+s+1},\infty}(L_2(\mathbb{R}^d))
\end{align*}
where $R_{j}:=D_{j}\Delta^{-\frac{1}{2}}$ is the Riesz transform on $\mathbb{R}^{d}$ for $j=1,\cdots,d$.
\end{theorem}

\subsection{Double operator integral}
In this subsection, we demand the notion of double operator integrals as specifically used
in \cite{HK2003, PS-crelle}. In what follows, $H$ is a separable Hilbert space. For self-adjoint operators $A_{1},A_{2}$ on $H$, let's formally define the operator $T_{\psi}^{A_{1},A_{2}}$ with $\psi\in L_{\infty}(\mathbb{R}^{2})$ by setting
$$T_{\psi}^{A_{1},A_{2}}(x)=\int_{\mathbb{R}^{2}}\psi(\lambda,\mu)dE^{1}(\lambda)xdE^{2}(\mu),\quad x\in B(H),$$
where $E^{1},E^{2}$ are the spectral measures of $A_{1},A_{2}$ separably. We have below result.
\begin{lemma}\label{double integral bd}
For $0<\beta<m$, define 
$$\psi_{\beta,m}(\lambda,\mu)
=\frac{\lambda^{\beta}-\mu^{\beta}}{\lambda^{m}-\mu^{m}}\lambda^{\frac{m-\beta}{2}}\mu^{\frac{m-\beta}{2}}
,\quad\lambda,\mu>0.$$ 
For any positive operator $B\in B(H)$, we have 
$$\|T_{\psi_{\beta,m}}^{B,B}\|_{B(H)\rightarrow B(H)}<\infty.$$
\end{lemma}
\begin{proof}
Let $\frac{\lambda}{\mu}=e^{t}$ with $t\in\mathbb{R}$. We have
\begin{align*}
\psi_{\beta,m}(\lambda,\mu)
=\frac{(\frac{\lambda}{\mu})^{\beta}-1}{(\frac{\lambda}{\mu})^m-1}(\frac{\lambda}{\mu})^{\frac{m-\beta}{2}}
=\frac{\sinh(\frac{\beta}{2}t)}{\sinh(\frac{m}{2}t)}
\end{align*}
where $t\mapsto\sinh(t)$ is the hyperbolic sine function. As $\beta<m$, the function $t\mapsto\frac{\sinh(\frac{\beta}{2}t)}{\sinh(\frac{m}{2}t)}$ belongs to the Schwartz class of $\mathbb{R}$. Let $g_{\beta,m}$ be its Fourier transform. By Fourier inversion,
\begin{align*}
\psi_{\beta,m}(\lambda,\mu)=\int_{\mathbb{R}}g_{\beta,m}(s)e^{2\pi\ri s\log_{e}(\frac{\lambda}{\mu})}ds
=\int_{\mathbb{R}}g_{\beta,m}(s)\lambda^{2\pi\ri s}\mu^{-2\pi\ri s}ds.
\end{align*}
Therefore, for $A\in B(H)$,
\begin{align*}
T_{\psi_{\beta,m}}^{B,B}(A)=\int_{\mathbb{R}}g_{\beta,m}(s)B^{2\pi\ri s}AB^{-2\pi\ri s}ds.
\end{align*}
Since $g_{\beta,m}$ is integrable in $L_{1}(\mathbb{R})$, the boundedness of $T_{\psi_{\beta,m}}^{B,B}$ follows.
\end{proof}

The following lemma is a direct result from functional calculus. Its a special case $m=2$ is established in \cite[lemma 5.4]{MSZ2023} so we omit the calculations.
\begin{lemma}\label{double integral expansion}
Let $B,A\in B(H)$ and assume that there is a constant $c>0$ such that $c\leq B\leq c^{-1}$. Let $\alpha,\gamma\geq0$, and $0<\beta<m$. Define $\psi_{s,m}$ as in Lemma \ref{double integral bd}. We have
$$B^{-\alpha}[A,B^{\beta}]B^{-\gamma}
=T_{\psi_{\beta,m}}^{B,B}(B^{\frac{\beta-m}{2}-\alpha}[A,B^{m}]B^{\frac{\beta-m}{2}-\gamma}).$$
\end{lemma}

\begin{corollary}\label{double integral co}
Let $A,B$ be positive operators on Hilbert space $H$ such that $A:{\rm dom}(B^{m})\rightarrow {\rm dom}(B^{m})$, $m>0$. Assume that $A$ is bounded on $H$ and ${\rm ker}(B)=\{0\}$. For $0<\beta<m$ and $\alpha,\gamma\in\mathbb{R}$, we have
$$\|B^{-\alpha}[A,B^{\beta}]B^{-\gamma}\|_{B(H)}\leq c_{\beta,m}
\|B^{\frac{\beta-m}{2}-\alpha}[A,B^{m}]B^{\frac{\beta-m}{2}-\gamma}\|_{B(H)}$$
with constant $c_{\beta,m}$ depending on $\beta,m$.
\end{corollary}
\begin{proof}
For each $n\geq1$, consider projection $P_{n}=\chi_{(\frac{1}{n},n)}(B)$ and Hilbert space $H_{n}=P_{n}(H)$. Let $B_{n}=B|_{H_{n}}$ and $A_{n}=P_{n}AP_{n}|_{H_{n}}$. Clearly, $A_{n}\in B(H_{n})$ and $B_{n}\in B(H_{n})$ satisfies the conditions of Lemma \ref{double integral expansion}. 
Hence,
$$B_{n}^{-\alpha}[A_{n},B_{n}^{\beta}]B_{n}^{-\gamma}
=T_{\psi_{\beta,m}}^{B_{n},B_{n}}(B_{n}^{\frac{\beta-m}{2}-\alpha}[A_{n},B_{n}^{m}]B_{n}^{\frac{\beta-m}{2}-\gamma}).$$
Using Lemma \ref{double integral bd},
$$\|B_{n}^{-\alpha}[A_{n},B_{n}^{\beta}]B_{n}^{-\gamma}\|_{B(H)}\leq c_{\beta,m}
\|B_{n}^{\frac{\beta-m}{2}-\alpha}[A_{n},B_{n}^{m}]B_{n}^{\frac{\beta-m}{2}-\gamma}\|_{B(H)}.$$
Observe that 
$$B_{n}^{-\alpha}[A_{n},B_{n}^{\beta}]B_{n}^{-\gamma}=
P_{n}\cdot B^{-\alpha}[A,B^{\beta}]B^{-\gamma}\cdot P_{n}$$
and
$$B_{n}^{\frac{\beta-m}{2}-\alpha}[A_{n},B_{n}^{m}]B_{n}^{\frac{\beta-m}{2}-\gamma}=
P_{n}\cdot B^{\frac{\beta-m}{2}-\alpha}[A,B^{m}]B^{\frac{\beta-m}{2}-\gamma}\cdot P_{n}.$$
Therefore,
$$\|P_{n}\cdot B^{-\alpha}[A,B^{\beta}]B^{-\gamma}\cdot P_{n}\|_{B(H)}\leq c_{\beta,m}
\|B^{\frac{\beta-m}{2}-\alpha}[A,B^{m}]B^{\frac{\beta-m}{2}-\gamma}\|_{B(H)}$$
for every $n\geq1$. The assertion follows from the Fatou property in $B(H)$.
\end{proof}

\subsection{Proof of Theorem \ref{euclidean compactness of commutators}}
Before proving Theorem \ref{euclidean compactness of commutators}, we establish the following necessary lemmas.
\begin{lemma}\label{pre-hadamard lemma} Let $P$ and $\phi$ be given as in Theorem \ref{euclidean compactness of commutators}. For every $n\in\mathbb{N},$ we have
$$(1+P)^{-\frac{2n-1}{2}}M_{\phi}D_k(1+P)^{n-1},(1+P)^nM_{\phi}D_k(1+P)^{-\frac{2n+1}{2}}\in B(L_2(\mathbb{R}^{d})).$$
This should be understood as follows: operators are well defined on $\bigcap_{k\geq0}W^{2k,2}(\mathbb{R}^d)$ and extend to bounded operators on $L_2(\mathbb{R}^d).$
\end{lemma}
\begin{proof} By the functional calculus,
$$(1+P)^{\alpha}:\bigcap_{k\geq0}{\rm dom}((1+P)^k)\to\bigcap_{k\geq0}{\rm dom}((1+P)^k),\quad z\in\mathbb{C}.$$
Hence, by Theorem \ref{standard elliptic theorem}, 
$$(1+P)^{\alpha}:\bigcap_{k\geq0}W^{2k,2}(\mathbb{R}^d)\to\bigcap_{k\geq0}W^{2k,2}(\mathbb{R}^d),\quad z\in\mathbb{C}.$$
Clearly, also 
$$M_{\phi}:\bigcap_{k\geq0}W^{2k,2}(\mathbb{R}^d)\to\bigcap_{k\geq0}W^{2k,2}(\mathbb{R}^d).$$
Hence, the operators in question are well defined on $\cap_{k\geq0}W^{2k,2}(\mathbb{R}^d).$
	
We write the second operator as
$$\Big((1+P)^nM_{\phi}D_k(1+P)^{-\frac{2n+1}{2}}:L_2(\mathbb{R}^d)\to L_2(\mathbb{R}^d)\Big)$$
$$=\Big((1+P)^nM_{\phi}D_k:W^{2n+1,2}(\mathbb{R}^d)\to L_2(\mathbb{R}^d)\Big)\circ$$
$$\circ \Big((1+P)^{-\frac{2n+1}{2}}:L_2(\mathbb{R}^d)\to W^{2n+1,2}(\mathbb{R}^d)\Big).$$
The first factor on the right hand side is obviously bounded. The second factor on the right hand side is bounded by Theorem \ref{standard elliptic theorem}. Hence, our second operator is bounded.

For every $\xi,\eta\in \cap_{k\geq0}W^{2k,2}(\mathbb{R}^d),$ we have
$$\langle (1+P)^{-\frac{2n-1}{2}}M_{\phi}D_k(1+P)^{n-1}\xi,\eta\rangle=\langle (1+P)^{n-1}D_kM_{\bar{\phi}} (1+P)^{-\frac{2n-1}{2}}\xi,\eta\rangle.$$
That the operator
$$(1+P)^{n-1}D_kM_{\bar{\phi}} (1+P)^{-\frac{2n-1}{2}}$$
is bounded can be seen as in the preceding paragraph. Thus, our first operator is also bounded.
\end{proof}

\begin{lemma}\label{from hadamard lemma} 
Let $P$ and $\phi$ be given as in Theorem \ref{euclidean compactness of commutators}. For every $z\in\mathbb{C},$ we have
$$(1+P)^{-\frac{z}{2}}M_{\phi}D_k(1+P)^{\frac{z-1}{2}},
(1+P)^{-\frac{z}{2}}M_{\phi}(1+P)^{\frac{z-1}{2}}\in B(L_{2}(\mathbb{R}^{d})).$$
This should be understood as follows: operators are well defined on $\bigcap_{k\geq0}W^{2k,2}(\mathbb{R}^d)$ and extend to bounded operators on $L_2(\mathbb{R}^d).$
\end{lemma}
\begin{proof} By the functional calculus,
$$(1+P)^z:\bigcap_{k\geq0}{\rm dom}((1+P)^k)\to\bigcap_{k\geq0}{\rm dom}((1+P)^k),\quad z\in\mathbb{C}.$$
Hence, by Theorem \ref{standard elliptic theorem}, 
$$(1+P)^z:\bigcap_{k\geq0}W^{2k,2}(\mathbb{R}^d)\to\bigcap_{k\geq0}W^{2k,2}(\mathbb{R}^d),\quad z\in\mathbb{C}.$$
Clearly, also 
$$M_{\phi}:\bigcap_{k\geq0}W^{2k,2}(\mathbb{R}^d)\to\bigcap_{k\geq0}W^{2k,2}(\mathbb{R}^d).$$
Hence, the operator in question is well defined on $\cap_{k\geq0}W^{2k,2}(\mathbb{R}^d).$
	
For $\xi,\eta\in \cap_{k\geq0}W^{2k,2}(\mathbb{R}^d),$ define entire function
$$F_{\xi,\eta}:z\to \langle (1+P)^{-\frac{z}{2}}M_{\phi}D_k(1+P)^{\frac{z-1}{2}}\xi,\eta\rangle,\quad z\in\mathbb{C}.$$

Fix $n\in\mathbb{N}$ and consider the function $F_{\xi,\eta}$ on the strip $\{-2n\leq\Re(z)\leq 2n-1\}.$ We have
$$F_{\xi,\eta}(z)=\langle M_{\phi}D_k(1+P)^{\frac{z-1}{2}}\xi,(1+P)^{-\frac{\bar{z}}{2}}\eta\rangle.$$
Hence,
\begin{align*}
&|F_{\xi,\eta}(z)|\\
&\leq\Big\|M_{\phi}D_k(1+P)^{\frac{z-1}{2}}\xi\Big\|_{L_2(\mathbb{R}^d)}
\Big\|(1+P)^{-\frac{\bar{z}}{2}}\eta\Big\|_{L_2(\mathbb{R}^d)}\\
&\leq \Big\|M_{\phi}D_k\Big\|_{W^{1,2}(\mathbb{R}^d)\to L_2(\mathbb{R}^d)}\Big\|(1+P)^{\frac{z-1}{2}}\xi\Big\|_{W^{1,2}(\mathbb{R}^d)}
\Big\|(1+P)^{-\frac{\bar{z}}{2}}\eta\Big\|_{L_2(\mathbb{R}^d)}\\
&\stackrel{Thm.\ref{standard elliptic theorem}}{\leq}c_P\Big\|M_{\phi}D_k\Big\|_{W^{1,2}(\mathbb{R}^d)\to L_2(\mathbb{R}^d)}\Big\|(1+P)^{\frac{z-1}{2}}\xi\Big\|_{{\rm dom}((P+1)^{\frac12})}\Big\|(1+P)^{-\frac{\bar{z}}{2}}\eta\Big\|_{L_2(\mathbb{R}^d)}\\
&\leq c_P\Big\|M_{\phi}D_k\Big\|_{W^{1,2}(\mathbb{R}^d)\to L_2(\mathbb{R}^d)}\|\xi\|_{{\rm dom}((P+1)^{\frac{2n-1}{2}})}\|\eta\|_{{\rm dom}(P^n)}\\	
&\stackrel{Thm.\ref{standard elliptic theorem}}{\leq} c_{P,n}\Big\|M_{\phi}D_k\Big\|_{W^{1,2}(\mathbb{R}^d)\to L_2(\mathbb{R}^d)}\|\xi\|_{W^{2n-1,2}(\mathbb{R}^d)}\|\eta\|_{W^{2n,2}(\mathbb{R}^d)}.
\end{align*}
Thus, $F_{\xi,\eta}$ is bounded on the strip $\{-2n\leq\Re(z)\leq 2n-1\}.$

Clearly,
$$F_{\xi,\eta}(2n-1+it)=\langle (1+P)^{-\frac{2n-1}{2}}M_{\phi}D_k(1+P)^{n-1}\big((1+P)^{\frac{it}{2}}\xi\big),(1+P)^{\frac{it}{2}}\eta\rangle.$$
We, therefore, have
\begin{align*}
&|F_{\xi,\eta}(2n-1+it)|\\
&\leq\Big\| (1+P)^{-\frac{2n-1}{2}}M_{\phi}D_k(1+P)^{n-1}\big((1+P)^{\frac{it}{2}}\xi\big)
\Big\|_{L_2(\mathbb{R}^d)}\Big\|(1+P)^{\frac{it}{2}}\eta\Big\|_{L_2(\mathbb{R}^d)}\\
&\leq\Big\| (1+P)^{-\frac{2n-1}{2}}M_{\phi}D_k(1+P)^{n-1}\Big\|_{L_2(\mathbb{R}^d)\to L_2(\mathbb{R}^d)}\times\\
&\qquad\times\Big\|(1+P)^{\frac{it}{2}}\xi\Big\|_{L_2(\mathbb{R}^d)}
\Big\|(1+P)^{\frac{it}{2}}\eta\Big\|_{L_2(\mathbb{R}^d)}\\
&\leq\Big\| (1+P)^{-\frac{2n-1}{2}}M_{\phi}D_k(1+P)^{n-1}\Big\|_{L_2(\mathbb{R}^d)\to L_2(\mathbb{R}^d)}\|\xi\|_{L_2(\mathbb{R}^d)}\|\eta\|_{L_2(\mathbb{R}^d)}.
\end{align*}
Here, the finiteness of the right hand side is guaranteed by Lemma \ref{pre-hadamard lemma}.

Similarly,
\begin{align*}
&|F_{\xi,\eta}(-2n+it)|\\
&\leq\Big\| (1+P)^nM_{\phi}D_k(1+P)^{-\frac{2n+1}{2}}\Big\|_{L_2(\mathbb{R}^d)\to L_2(\mathbb{R}^d)}\|\xi\|_{L_2(\mathbb{R}^d)}\|\eta\|_{L_2(\mathbb{R}^d)}.
\end{align*}

By the Maximum Modulus Principle,
\begin{align*}
&|F_{\xi,\eta}(z)|\\
&\leq\Big(\Big\| (1+P)^{-\frac{2n-1}{2}}M_{\phi}D_k(1+P)^{n-1}\Big\|_{L_2(\mathbb{R}^d)\to L_2(\mathbb{R}^d)}+\\
&+\Big\| (1+P)^nM_{\phi}D_k(1+P)^{-\frac{2n+1}{2}}\Big\|_{L_2(\mathbb{R}^d)\to L_2(\mathbb{R}^d)}\Big)\times \|\xi\|_{L_2(\mathbb{R}^d)}\|\eta\|_{L_2(\mathbb{R}^d)}
\end{align*}
for every $z$ in the strip $\{-2n\leq\Re(z)\leq 2n-1\}.$ Hence,
\begin{align*}
&\Big\|(1+P)^{-\frac{z}{2}}M_{\phi}D_k(1+P)^{\frac{z-1}{2}}\Big\|_{L_2(\mathbb{R}^d)\to L_2(\mathbb{R}^d)}\\
&\leq\Big\| (1+P)^{-\frac{2n-1}{2}}M_{\phi}D_k(1+P)^{n-1}\Big\|_{L_2(\mathbb{R}^d)\to L_2(\mathbb{R}^d)}+\\
&+\Big\| (1+P)^nM_{\phi}D_k(1+P)^{-\frac{2n+1}{2}}\Big\|_{L_2(\mathbb{R}^d)\to L_2(\mathbb{R}^d)}
\end{align*}
for every $z$ in the strip $\{-2n\leq\Re(z)\leq 2n-1\}.$ This completes the proof.
\end{proof}

\begin{lemma}\label{P commutator bd-integer}
Let $P$ and $\phi$ be given as in Theorem \ref{euclidean compactness of commutators}.  If $\alpha,\gamma\in\mathbb{R}$ satisfy $\alpha+\gamma=1$, then
$$(1+P)^{-\frac{\alpha}{2}}[M_{\phi},P](1+P)^{-\frac{\gamma}{2}}
\in B(L_{2}(\mathbb{R}^{d})).$$
\end{lemma}
\begin{proof} 
Write $$P=\sum_{k,l=1}^dM_{a_{kl}}D_kD_l+\sum_{k=1}^dM_{b_k}D_k+M_f,$$
where functions $\{a_{kl}\}_{k,l=1}^d,$ $\{b_k\}_{k=1}^d$ and $f$ are smooth and constant outside some ball. We have
\begin{align*}
[M_{\phi},P]
&=\sum_{k,l=1}^dM_{a_{kl}}[M_{\phi},D_kD_l]+\sum_{k=1}^dM_{b_k}[M_{\phi},D_k]\\
&=-\sum_{k,l=1}^dM_{a_{kl}}\cdot M_{D_k\phi}D_l-\sum_{k,l=1}^dM_{a_{kl}}\cdot D_kM_{D_l\phi}
-\sum_{k=1}^dM_{b_k}\cdot M_{D_k\phi}\\
&=-\sum_{k,l=1}^dM_{a_{kl}}\cdot M_{D_k\phi}D_l-\sum_{k,l=1}^dM_{a_{kl}}\cdot M_{D_l\phi}D_k-\sum_{k,l=1}^dM_{a_{kl}}\cdot M_{D_kD_l\phi}+\\
&\qquad\qquad\qquad\qquad\qquad\qquad\qquad\qquad\qquad\qquad\qquad
-\sum_{k=1}^dM_{b_k}\cdot M_{D_k\phi}.
\end{align*}
The assertion follows now from Lemma \ref{from hadamard lemma}.
\end{proof}

\begin{lemma}\label{tab lemma} Let $\alpha,\beta\in\mathbb{R}.$ Set
$$T_{\alpha,\beta}=(1+P)^{-\frac{\alpha}{2}}[M_{\phi},(1+P)^{\frac{\beta}{2}}](1+P)^{\frac{\alpha-\beta+1}{2}}.$$
We have
$$T_{\alpha,\beta+1}=T_{\alpha,\beta}+T_{\alpha-\beta,1}.$$
\end{lemma}
\begin{proof} By the Leibniz rule,
\begin{align*}
&T_{\alpha,\beta+1}
=(1+P)^{-\frac{\alpha}{2}}[M_{\phi},(1+P)^{\frac{\beta}{2}}\cdot (1+P)^{\frac12}](1+P)^{\frac{\alpha-\beta}{2}}\\
&=(1+P)^{-\frac{\alpha}{2}}[M_{\phi},(1+P)^{\frac{\beta}{2}}](1+P)^{\frac{\alpha-\beta+1}{2}}
+(1+P)^{\frac{\beta-\alpha}{2}}[M_{\phi}, (1+P)^{\frac12}](1+P)^{\frac{\alpha-\beta}{2}}\\
&=T_{\alpha,\beta}+T_{\alpha-\beta,1}.
\end{align*}
The assertion follows.
\end{proof}

\begin{lemma}\label{P-Mf boundedess}
Let $P$ and $\phi$ be given as in Theorem \ref{euclidean compactness of commutators}. If $\alpha,\beta,\gamma\in\mathbb{R}$ satisfy $\alpha+\gamma+1=\beta$, then
$$(1+P)^{-\frac{\alpha}{2}}[M_{\phi},(1+P)^{\frac{\beta}{2}}](1+P)^{-\frac{\gamma}{2}}
\in B(L_{2}(\mathbb{R}^{d})).$$
\end{lemma}
\begin{proof} For $\beta=0,$ there is nothing to prove.
	
{\bf Case 1:} Let $0\leq \beta<2.$ For $\beta=0,$ there is nothing to prove and, so, we assume that $0<\beta<2.$ Set $A=M_{\phi}$ and $B=(1+P)^{\frac{1}{2}}.$ By Corollary \ref{double integral co} and Lemma \ref{P commutator bd-integer},
\begin{align*}
&\|(1+P)^{-\frac{\alpha}{2}}
[M_{\phi},(1+P)^{\frac{\beta}{2}}](1+P)^{-\frac{\gamma}{2}}\|_{B(L_{2}(\mathbb{R}^{d}))}\\
&\leq c_{\beta}\|(1+P)^{\frac{\beta-2}{4}-\frac{\alpha}{2}}[M_{\phi},P]
(1+P)^{\frac{\beta-2}{4}-\frac{\gamma}{2}}\|_{B(L_{2}(\mathbb{R}^{d}))}
<\infty.
\end{align*}
This proves the assertion for $0<\beta<2.$

{\bf Case 2:} Consider now the general case. By Lemma \ref{tab lemma}, we have 
$$T_{\alpha,\beta+1}=T_{\alpha,\beta}+T_{\alpha-\beta,1},\quad T_{\alpha,\beta-1}=T_{\alpha,\beta}-T_{\alpha-\beta+1,1}.$$
By induction, we can express $T_{\alpha,\beta}$ as linear combination of $T_{\delta,\gamma}$'s with $\gamma\in[0,1].$ The assertion in general case follows now from that in Case 1.
\end{proof}

\begin{lemma}\label{P-Mf compactness}
Let $P$ and $\phi$ be given as in Theorem \ref{euclidean compactness of commutators}. If $\alpha,\gamma\geq0$ with $\alpha+\gamma>0$, then
$$(1+P)^{-\frac{\alpha}{2}}M_{\phi}(1+P)^{-\frac{\gamma}{2}}
\in\mathcal{L}_{\frac{d}{\alpha+\gamma},\infty}(L_{2}(\mathbb{R}^{d})).$$
\end{lemma}
\begin{proof}
Write 
\begin{align*}
(1+P)^{-\frac{\alpha}{2}}M_{\phi}(1+P)^{-\frac{\gamma}{2}}
=(1+P)^{-\frac{\alpha}{2}}(1+\Delta)^{\frac{\alpha}{2}}\cdot
(1-\Delta)^{-\frac{\alpha}{2}}M_{\phi}(1+\Delta)^{-\frac{\gamma}{2}}\cdot\\
\cdot(1+\Delta)^{\frac{\gamma}{2}}(1+P)^{-\frac{\gamma}{2}}.
\end{align*}
The first and the third factors on the right hand side are bounded by Theorem \ref{standard elliptic theorem}. The second factor on the right hand side belongs to $\mathcal{L}_{\frac{d}{\alpha+\gamma},\infty}(L_{2}(\mathbb{R}^{d}))$ by Lemma \ref{symmetric Cwikel estimate}.
\end{proof}

It is time to prove Theorem \ref{euclidean compactness of commutators}.
\begin{proof}[Proof of Theorem \ref{euclidean compactness of commutators}]
Write $\phi=\phi_{1}\phi_{2}$ with $\phi_{1},\phi_{2}\in C_{c}^{\infty}(\mathbb{R}^{d})$. We have
\begin{align*}
&(1+P)^{-\frac{\alpha}{2}}[M_{\phi},(1+P)^{-\frac{\beta}{2}}](1+P)^{-\frac{\gamma}{2}}\\
&=(1+P)^{-\frac{\alpha}{2}}M_{\phi_{1}}(1+P)^{-\frac{\beta+\gamma+1}{2}}
\cdot(1+P)^{\frac{\beta+\gamma+1}{2}}[M_{\phi_{2}},(1+P)^{-\frac{\beta}{2}}](1+P)^{-\frac{\gamma}{2}}+\\
&+(1+P)^{-\frac{\alpha}{2}}[M_{\phi_{1}},(1+P)^{-\frac{\beta}{2}}](1+P)^{\frac{\beta+\alpha+1}{2}}
\cdot(1+P)^{-\frac{\beta+\alpha+1}{2}}M_{\phi_{2}}(1+P)^{-\frac{\gamma}{2}}.
\end{align*}
In the first summand on the right hand side, the first factor belongs to the ideal $\mathcal{L}_{\frac{1}{\alpha+\beta+\gamma+1},\infty}(L_{2}(\mathbb{R}^{d}))$ by Lemma \ref{P-Mf compactness}, while the second factor is bounded by Lemma \ref{P-Mf boundedess}. In the second summand on the right hand side, the first factor is bounded by Lemma \ref{P-Mf boundedess}, while the second factor belongs to $\mathcal{L}_{\frac{1}{\alpha+\beta+\gamma+1},\infty}(L_{2}(\mathbb{R}^{d}))$ by Lemma \ref{P-Mf compactness}. This completes the proof.
\end{proof}

\subsection{Proof of Theorem \ref{1 over 2 case}}
In this subsection, we prove Theorem \ref{1 over 2 case}. 
\begin{lemma}\label{bd in second sum}
Let $P$ and $\phi$ be given in Theorem \ref{1 over 2 case}. Given a smooth function $f$ on $\mathbb{R}^{d}$ and $0\leq s\leq2$, we have
$$(1+P)^{-\frac{s}{4}}[P,M_{\phi}R_{j}M_{\phi}](1+P)^{\frac{s}{4}-\frac{1}{2}}
\in B(L_{2}(\mathbb{R}^{d})).$$
\end{lemma}
\begin{proof}
Denote 
$$A_{1}=M_{\phi}D_{j}(1+\Delta)^{-\frac{1}{2}}M_{\phi}\quad\mbox{and}\quad
A_{2}=M_{\phi}D_{j}(\Delta^{-\frac{1}{2}}-(1+\Delta)^{-\frac{1}{2}})M_{\phi}.$$
Write
\begin{align*}
[P,M_{\phi}R_{j}M_{\phi}]=[P,A_{1}]+[P,A_{2}].
\end{align*}
By the assumption on $P$, we have $P\in\Psi^{2}(\mathbb{R}^{d})$.
Since $D_{j}(1+\Delta)^{-\frac{1}{2}}\in\Psi^{0}(\mathbb{R}^{d})$ and $\phi\in C_{c}^{\infty}(\mathbb{R}^{d})$, it follows from \eqref{pseudo-differential op} that
$$(1+\Delta)^{-\frac{s}{4}}[P,A_{1}](1+\Delta)^{\frac{s}{4}-\frac{1}{2}}\in\Psi^{0}(\mathbb{R}^{d}),$$
where we use the fact that
$$(1+\Delta)^{\frac{t}{2}}\in\Psi^{t}(\mathbb{R}^{d}),\quad\forall\; t\in\mathbb{R}.$$
By \cite[Theorem 1,page 235]{Stein-Murphy1993}, $(1+\Delta)^{-\frac{s}{4}}[P,A_{1}](1+\Delta)^{\frac{s}{4}-\frac{1}{2}}$ is bounded on $L_{2}(\mathbb{R}^{d})$. Therefore, the operator
\begin{align*}
&(1+P)^{-\frac{s}{4}}[P,A_{1}](1+P)^{\frac{s}{4}-\frac{1}{2}}\\
&=(1+P)^{-\frac{s}{4}}(1+\Delta)^{\frac{s}{4}}
\cdot(1+\Delta)^{-\frac{s}{4}}[P,A_{1}](1+\Delta)^{\frac{s}{4}-\frac{1}{2}}\cdot
(1+\Delta)^{\frac{1}{2}-\frac{s}{4}}(1+P)^{\frac{s}{4}-\frac{1}{2}}
\end{align*}
is bounded on $L_{2}(\mathbb{R}^{d})$ because the boundedness for the first and third factors on the right hand side of the equality follows from Theorem \ref{standard elliptic theorem}, for the second one follows by previous argument. Turn to $[P,A_{2}]$. 
Observe that
\begin{align*}
&(1+P)^{-\frac{s}{4}}PA_{2}(1+P)^{\frac{s}{4}-\frac{1}{2}}\\
&=(1+P)^{-\frac{s}{4}}(1+\Delta)^{\frac{s}{4}}\cdot(1+\Delta)^{-\frac{s}{4}}
PM_{\phi}(1+\Delta)^{\frac{s}{4}-\frac{3}{2}}
\cdot\frac{R_{j}\sqrt{1+\Delta}}{\sqrt{1+\Delta}+\sqrt{\Delta}}\cdot\\
&\quad\cdot(1+\Delta)^{\frac{1}{2}-\frac{s}{4}}M_{\phi}(1+\Delta)^{\frac{s}{4}-\frac{1}{2}}\cdot
(1+\Delta)^{\frac{1}{2}-\frac{s}{4}}(1+P)^{\frac{s}{4}-\frac{1}{2}}.
\end{align*}
On right hand side of the equality, the first and fifth factors are bounded by Theorem \ref{standard elliptic theorem}, the second and fourth factors belong to $\Psi^{0}(\mathbb{R}^{d})$ and therefore are bounded, the third factor is obviously bounded. Thus, $(1+P)^{-\frac{s}{4}}PA_{2}(1+P)^{\frac{s}{4}-\frac{1}{2}}$ is a bounded operator. 
Similarly, 
$(1+P)^{-\frac{s}{4}}A_{2}P(1+P)^{\frac{s}{4}-\frac{1}{2}}$ is also bounded.
In a sum, the operator $(1+P)^{-\frac{s}{4}}[P,A_{2}](1+P)^{\frac{s}{4}-\frac{1}{2}}$ is bounded on $L_{2}(\mathbb{R}^{d})$. The proof is complete.
\end{proof}

\begin{lemma}\label{boundedness of Mf-Riesz}
Let $P$ and $\phi$ be given in Theorem \ref{1 over 2 case}. If $s\in(0,2)$, then
\begin{align*}
(1+P)^{\frac{1-s}{2}}[(1+P)^{\frac{s}{2}},M_{\phi}R_{j}M_{\phi}]\in B(L_{2}(\mathbb{R}^{d})).
\end{align*}
\end{lemma}
\begin{proof}
Using Corollary \ref{double integral co}, there is a constant $c_{s}$ relying on $s$ such that  $$\|(1+P)^{\frac{1-s}{2}}[(1+P)^{\frac{s}{2}},M_{\phi}R_{j}M_{\phi}]\|_{\infty}\leq c_{s}\|(1+P)^{-\frac{s}{4}}[P ,M_{\phi}R_{j}M_{\phi}](1+P)^{\frac{s}{4}-\frac{1}{2}}\|_{\infty}.$$
By Lemma \ref{bd in second sum}, the second factor on the inequality is finite. The proof is complete.
\end{proof}

Let us start the proof of Theorem \ref{1 over 2 case}.
\begin{proof}[Proof of Theorem \ref{1 over 2 case}]
Take $\psi\in C_{c}^{\infty}(\mathbb{R}^{d})$ such that $\phi=\phi\psi$ and therefore
$M_{\phi}R_{j}M_{\phi}=M_{\psi}\cdot M_{\phi}R_{j}M_{\phi}\cdot M_{\psi}.$
We have
\begin{align*}
&(1+P)^{-\frac{\alpha}{2}}[M_{\phi}R_{j}M_{\phi},(1+P)^{-\frac{s}{2}}](1+P)^{-\frac{\gamma}{2}}\\
&=(1+P)^{-\frac{\alpha}{2}}[M_{\psi},(1+P)^{-\frac{s}{2}}]\cdot M_{\phi}R_{j}\cdot M_{\phi}(1+P)^{-\frac{\gamma}{2}}+\\
&+(1+P)^{-\frac{\alpha}{2}}M_{\psi}[M_{\phi}R_{j}M_{\phi},(1+P)^{-\frac{s}{2}}]M_{\psi}(1+P)^{-\frac{\gamma}{2}}+\\
&+(1+P)^{-\frac{\alpha}{2}}M_{\phi}\cdot R_{j}M_{\phi}\cdot[M_{\psi},(1+P)^{-\frac{s}{2}}](1+P)^{-\frac{\gamma}{2}}.
\end{align*}
By Theorem \ref{euclidean compactness of commutators}, the first factor belongs to $\mathcal{L}_{\frac{d}{\alpha+s+1},\infty}\cdot\mathcal{L}_{\infty}
\cdot\mathcal{L}_{\frac{d}{\gamma},\infty}$ and therefore
$\mathcal{L}_{\frac{d}{\alpha+\gamma+s+1},\infty}$. By Theorem \ref{euclidean compactness of commutators}, third factor belong to $\mathcal{L}_{\frac{d}{\alpha},\infty}\cdot\mathcal{L}_{\infty}\cdot
\mathcal{L}_{\frac{d}{\gamma+s+1},\infty}$ and therefore
$\mathcal{L}_{\frac{d}{\alpha+\gamma+s+1},\infty}$
Here $\mathcal{L}_{\frac{d}{\alpha},\infty}$ (or $\mathcal{L}_{\frac{d}{\gamma},\infty}$) is understood as $\mathcal{L}_{\infty}$ when $\alpha$ (or $\gamma$) is $0$.
For the second factor, Lemma \ref{boundedness of Mf-Riesz} and Lemma \ref{P-Mf compactness} give that
\begin{align*}
&(1+P)^{-\frac{\alpha}{2}}M_{\psi}[M_{\phi}R_{j}M_{\phi},(1+P)^{-\frac{s}{2}}]
M_{\psi}(1+P)^{-\frac{\gamma}{2}}\\
&=(1+P)^{-\frac{\alpha}{2}}M_{\psi}(1+P)^{-\frac{1}{2}}\cdot
(1+P)^{\frac{1-s}{2}}[(1+P)^{\frac{s}{2}},M_{\phi}R_{j}M_{\phi}]\cdot\\
&\quad\cdot(1+P)^{-\frac{s}{2}}M_{\psi}(1+P)^{-\frac{\gamma}{2}}\\
&\in\mathcal{L}_{\frac{d}{\alpha+1},\infty}\cdot\mathcal{L}_{\infty}\cdot\mathcal{L}_{\frac{d}{\gamma+s},\infty}
\subset\mathcal{L}_{\frac{d}{\alpha+\gamma+s+1},\infty}.
\end{align*}
The proof is complete.
\end{proof}

\section{Asymptotics on the curved plane}\label{sec-4}
In this section, we generalize Theorem \ref{asymptotic on Rd} to Theorem \ref{for curved plane thm} by replacing the Laplacian with an elliptic differential operator on $\mathbb{R}^{d}$ for the index $p=\frac{d}{2}$.
\begin{theorem}\label{for curved plane thm} 
Let $P:W^{2,2}(\mathbb{R}^d)\to L_2(\mathbb{R}^d)$ be a positive, self-adjoint, elliptic differential operator whose coefficients are constant outside some ball. If $S\in\Pi_{\mathbb{R}^d}$ is compactly supported, then
$$\lim_{t\to\infty}t^{\frac2d}\mu(t,S(1+P)^{-1})=c_d\|{\rm sym}(S)\|_{L_{\frac{d}{2}}(\mathbb{R}^d\times\mathbb{R}^d,e^{-\sigma_P}dm)}$$
where $\sigma_P$ is the principal symbol of $P$.
\end{theorem}

Our starting point for proving Theorem \ref{for curved plane thm} is the following lemma.
\begin{lemma}\label{approximate lemma section2} Let $P:W^{2,2}(\mathbb{R}^d)\to L_2(\mathbb{R}^d)$ be a positive, self-adjoint, elliptic differential operator whose coefficients are constant outside some ball. There exists $T\in\Pi_{\mathbb{R}^d}$ such that, for every $\phi\in C^{\infty}_c(\mathbb{R}^d),$
$$(1+P)^{-1}M_{\phi}\in T(1+\Delta)^{-1}M_{\phi}+\mathcal{L}_{\frac{d}{3},\infty}(L_2(\mathbb{R}^d))$$
and
$${\rm sym}(T)(t,s)=(\sum_{i,j=1}^da_{ij}(t)s_is_j)^{-1},\quad t\in\mathbb{R}^d,\quad s\in\mathbb{S}^{d-1}$$
where $\{a_{ij}\}_{i,j=1}^{d}$ are the coefficients of $P$ with respect to derivatives $\{D_{i}D_{j}\}_{i,j=1}^{d}$.   
\end{lemma}
\begin{proof} Set 
$$\sigma(t,s)=\sum_{i,j=1}^da_{ij}(t)s_is_j,\quad t,s\in\mathbb{R}^d.$$
By the Complex Powers Theorem (see e.g. Theorem 1.6.22 in \cite{LMSZvol2}), we have
$$(1+P)^{-1}={\rm Op}((1+\sigma)^{-1})+P_1,\quad P_1\in \Psi^{-3}(\mathbb{R}^d).$$

Let $\psi\in C^{\infty}_c(\mathbb{R}^d)$ be such that $\psi=1$ near $0.$ Set
$$a(t,s)=(\sum_{i,j=1}^da_{ij}(t)\frac{s_i}{|s|}\frac{s_j}{|s|})^{-1},\quad t,s\in\mathbb{R}^d,$$
and
$$b(t,s)=a(t,s)\cdot (1-\psi)(s)\cdot (1+|s|^2)^{-1},\quad t,s\in\mathbb{R}^d.$$
Clearly,
$$(1+\sigma)^{-1}-b\in S^{-3}(\mathbb{R}^d\times\mathbb{R}^d).$$
Thus,
$${\rm Op}((1+\sigma)^{-1})-{\rm Op}(b)\in \Psi^{-3}(\mathbb{R}^d).$$
Using the preceding paragraph, we write
$$(1+P)^{-1}={\rm Op}(b)+P_2,\quad P_2\in \Psi^{-3}(\mathbb{R}^d).$$
Next,
$${\rm Op}(b)={\rm Op}(a)\cdot (1-\psi)(\nabla)\cdot (1+\Delta)^{-1}.$$
Thus,
$$(1+P)^{-1}={\rm Op}(a)\cdot (1-\psi)(\nabla)\cdot (1+\Delta)^{-1}+P_2,\quad P_2\in \Psi^{-3}(\mathbb{R}^d).$$
Clearly,
$$P_2M_{\phi}\in\mathcal{L}_{\frac{d}{3},\infty}(L_2(\mathbb{R}^d)).$$
Hence,
$$(1+P)^{-1}M_{\phi}\in {\rm Op}(a)\cdot (1-\psi)(\nabla)\cdot (1+\Delta)^{-1}M_{\phi}+\mathcal{L}_{\frac{d}{3},\infty}(L_2(\mathbb{R}^d)).$$
It remains to note that ${\rm Op}(a)\in\Pi_{\mathbb{R}^d}$ and its principal symbol
$${\rm sym}({\rm Op}(a))(t,s)=(\sum_{i,j=1}^da_{ij}(t)s_is_j)^{-1},\quad t\in\mathbb{R}^d,\quad s\in\mathbb{S}^{d-1}.$$
In addition,
$$(1-\psi)(\nabla)\cdot (1+\Delta)^{-1}M_{\phi}\in (1+\Delta)^{-1}M_{\phi}+\mathcal{L}_{\frac{d}{3},\infty}(L_2(\mathbb{R}^d)).$$
This completes the proof.
\end{proof}

\begin{lemma}\label{polar coordinates lemma} Let
$$\sigma(t,s)=\sum_{i,j=1}^da_{ij}(t)s_is_j,\quad t,s\in\mathbb{R}^d.$$
For every $F\in L_{\infty}(\mathbb{R}^d\times\mathbb{S}^{d-1}),$ we have 
\footnote{Recall that we identify functions on the sphere with homogeneous functions on Euclidean space.}
$$\|F\cdot \sigma^{-1}\|_{L_{\frac{d}{2}}(\mathbb{R}^d\times\mathbb{S}^{d-1})}
=\big(\frac12\Gamma(\frac{d}{2})\big)^{-\frac2d}
\|F\|_{L_{\frac{d}{2}}(\mathbb{R}^d\times\mathbb{R}^d,e^{-\sigma}dm)}$$
where $\Gamma(p):=\int_{0}^{\infty}x^{p-1}e^{-p}dx$, $p>0$, is the Gamma function.
\end{lemma}
\begin{proof} We have
$$\|F\|_{L_{\frac{d}{2}}(\mathbb{R}^d\times\mathbb{R}^d,e^{-\sigma}dm)}^{\frac{d}{2}}=\int_{\mathbb{R}^d\times\mathbb{R}^d}|F|^{\frac{d}{2}}(t,s)e^{-\sigma(t,s)}dtds.$$
Let us pass to spherical coordinates $r=|s|,$ $u=\frac{s}{|s|}.$ We have
$$|F|^{\frac{d}{2}}(t,s)=|F|^{\frac{d}{2}}(t,u),\quad e^{-\sigma(t,s)}=e^{-r^2\sigma(t,u)},\quad ds=r^{d-1}dr du.$$
Thus,
$$\|F\|_{L_{\frac{d}{2}}(\mathbb{R}^d\times\mathbb{R}^d,e^{-\sigma}dm)}^{\frac{d}{2}}
=\int_{\mathbb{R}^d\times\mathbb{S}^{d-1}}|F|^{\frac{d}{2}}(t,u)
\Big(\int_0^{\infty}e^{-r^2\sigma(t,u)}r^{d-1}dr\Big)dtdu.$$
Clearly,
$$\int_0^{\infty}e^{-r^2\sigma(t,u)}r^{d-1}dr=
\sigma(t,u)^{-\frac{d}{2}}\cdot\int_0^{\infty}e^{-x^2}x^{d-1}dx
=\sigma(t,u)^{-\frac{d}{2}}\cdot\frac12\Gamma(\frac{d}{2}).$$
Finally,
$$\|F\|_{L_{\frac{d}{2}}(\mathbb{R}^d\times\mathbb{R}^d,e^{-\sigma}dm)}^{\frac{d}{2}}=\frac12\Gamma(\frac{d}{2})\cdot \int_{\mathbb{R}^d\times\mathbb{S}^{d-1}}|F|^{\frac{d}{2}}(t,u)\sigma(t,u)^{-\frac{d}{2}}dtdu.$$
This completes the proof.
\end{proof}

Now, we prove Theorem \ref{for curved plane thm}.
\begin{proof}[Proof of Theorem \ref{for curved plane thm}] 
Observe that $$\lim_{t\to\infty}t^{\frac2d}\mu(t,S(1+P)^{-1})=\lim_{t\to\infty}t^{\frac2d}\mu(t,|S|(1+P)^{-1}),\quad
S\in\Pi_{\mathbb{R}^{d}}.$$
Without loss of generality, Suppose that $S\geq0$.

Let $\phi\in C^{\infty}_c(\mathbb{R}^d)$ be such that $S=M_{\phi}S.$ Let $T\in\Pi_{\mathbb{R}^d}$ be given by Lemma \ref{approximate lemma section2}. 
	
By Lemma \ref{approximate lemma section2}, we have
$$(1+P)^{-1}S=(1+P)^{-1}M_{\phi}S
\in T(1+\Delta)^{-1}M_{\phi}S+\mathcal{L}_{\frac{d}{3},\infty}(L_2(\mathbb{R}^d)).$$	
By \cite[Theorem A.1]{FSZ2024endpoint},
$$[(1+\Delta)^{-1},M_{\phi}S]\in(\mathcal{L}_{\frac{d}{2},\infty})_0(L_2(\mathbb{R}^d)).$$
And, $\mathcal{L}_{\frac{d}{3},\infty}(L_2(\mathbb{R}^d))
\subset(\mathcal{L}_{\frac{d}{2},\infty})_0(L_2(\mathbb{R}^d))$. Hence,
$$(1+P)^{-1}S\in TM_{\phi}S(1+\Delta)^{-1}+(\mathcal{L}_{\frac{d}{2},\infty})_0(L_2(\mathbb{R}^d)).$$	
Clearly, $[T,M_{\phi}]$ is compact. Hence,
$$(1+P)^{-1}S\in M_{\phi}TS(1+\Delta)^{-1}+(\mathcal{L}_{\frac{d}{2},\infty})_0(L_2(\mathbb{R}^d)).$$	
The operator $M_{\phi}TS$ is compactly supported. By Theorem \ref{asymptotic on Rd}, there exists a limit
$$\lim_{t\rightarrow\infty} t^{\frac2d}\mu\Big(t,M_{\phi}TS(1+\Delta)^{-1}\Big)
=c_d\|\sym(M_{\phi}TS)\|_{L_{\frac{d}{2}}(\mathbb{R}^{d}\times\mathbb{S}^{d-1})}.$$
Clearly,
\begin{align*}
\sym(STM_{\phi})&=\sym(S)\sym(T)\sym(M_{\phi})\\
&=\sym(S)\sym(M_{\phi})\sym(T)\\
&=\sym(SM_{\phi})\sym(T)\\
&=\sym(S)\sym(T).
\end{align*}
Hence,
$$\lim_{t\rightarrow\infty} t^{\frac2d}\mu\Big(t,M_{\phi}TS(1+\Delta)^{-1}\Big)
=c_d\|\sym(S)\sym(T)\|_{L_{\frac{d}{2}}(\mathbb{R}^{d}\times\mathbb{S}^{d-1})}.$$
The assertion follows now from Lemma \ref{polar coordinates lemma}.
\end{proof}

\section{Powers of nearly commuting operators}\label{sec-5}
In this section, we prove Theorem \ref{main example thm}, the spectral characterization for powers of two nearly commutating operators, which substitutes the Complex Powers Theorem for pseudo-differential operators on manifold.
\begin{theorem}\label{main example thm} 
Let $(X,G)$ be a compact Riemannian manifold and let $\Delta_G$ be the Laplace-Beltrami operator. For every $0\leq S\in\Pi_X,$ we have
$$S^s(1+\Delta_G)^{-s}-(S^{\frac12}(1+\Delta_G)^{-1}S^{\frac12})^s
\in(\mathcal{L}_{\frac{d}{2s},\infty})_{0}(L_{2}(X,{\rm vol}_{G})),\quad s>0.$$
\end{theorem}

\subsection{Abstract estimate}
In this subsection, we prove Theorem \ref{nearly commuting powers}. Firstly, fix $H$ a separable Hilbert space throughout this subsection.
\begin{theorem}\label{nearly commuting powers} Let $0\leq A\in\mathcal{L}_{\infty}(H)$ and let $0\leq B\in\mathcal{L}_{q,\infty}(H)$ with $q>0.$ If $$[A^{\epsilon},B^{\delta}]\in(\mathcal{L}_{\frac{q}{\delta},\infty})_0(H)\quad\mbox{for any}\quad \epsilon,\delta>0,$$
then
$$B^sA^s-(A^{\frac12}BA^{\frac12})^s\in(\mathcal{L}_{\frac{q}{s},\infty})_0(H),\quad\mbox{if}\quad s>0.$$
\end{theorem}

The following Lemma is a special case of Theorem \ref{nearly commuting powers} established in \cite{CSZ} (see Lemma 5.3 there).

\begin{lemma}\label{CSZ main lemma} 
Let $0\leq A\in\mathcal{L}_{\infty}(H)$ and let $0\leq B\in\mathcal{L}_{p,\infty}(H)$ with $p>1.$ If $[A^{\frac12},B]\in(\mathcal{L}_{p,\infty})_0(H),$ then
$$B^pA^p-(A^{\frac12}BA^{\frac12})^p\in(\mathcal{L}_{1,\infty})_0(H).$$
\end{lemma}

\begin{lemma}\label{bks lemma} If $0\leq X,Y\in\mathcal{L}_{r,\infty}(H)$ and $X-Y\in(\mathcal{L}_{r,\infty})_0(H)$ with $r>0$, then
$$X^t-Y^t\in(\mathcal{L}_{\frac{r}{t},\infty})_0(H),\quad t>0.$$
\end{lemma}
\begin{proof} 
Fix $n\in\mathbb{N}$ such that $n>t.$ Since $X-Y\in(\mathcal{L}_{r,\infty})_0(H),$ applying the proof of \cite[Corollary 1.3.15]{LMSZvol2} and the most general form of the Birman-Koplienko-Solomyak inequality established in \cite{HSZ}, we write
$$X^{\frac{t}{n}}-Y^{\frac{t}{n}}\in(\mathcal{L}_{\frac{nr}{t},\infty})_0(H).$$
Now, we write
$$X^t-Y^t=\sum_{k=0}^{n-1}X^{\frac{kt}{n}}(X^{\frac{t}{n}}-Y^{\frac{t}{n}})Y^{\frac{(n-1-k)t}{n}}.$$
The assertion follows now from H\"older inequality.
\end{proof}

The next lemma is an easy application of Lemma \ref{CSZ main lemma} and Lemma \ref{bks lemma}.

\begin{proof}[Proof of Theorem \ref{nearly commuting powers}] Choose $n\in\mathbb{N}$ such that $ns>1.$ Choose $m\in\mathbb{N}$ such that $mq>ns$ and $mq>1.$

By assumption, $0\leq A^{\frac{ns}{m}}\in\mathcal{L}_{\infty}(H)$ and $0\leq B^{\frac{ns}{m}}\in\mathcal{L}_{\frac{mq}{ns},\infty}(H)$ are such that $$[A^{\frac{ns}{2m}},B^{\frac{ns}{m}}]\in(\mathcal{L}_{\frac{mq}{ns},\infty})_0(H).$$
Applying Lemma \ref{CSZ main lemma} with $p=\frac{mq}{ns}$ to the couple $(A^{\frac{ns}{m}},B^{\frac{ns}{m}}),$ we obtain
\begin{equation}\label{postdenis eq0}
B^qA^q-(A^{\frac{ns}{2m}}B^{\frac{ns}{m}}A^{\frac{ns}{2m}})^{\frac{mq}{ns}}
\in(\mathcal{L}_{1,\infty})_0(H).
\end{equation}

By assumption, $0\leq A^{\frac{n}{m}}\in\mathcal{L}_{\infty}(H)$ and $0\leq B^{\frac{n}{m}}\in\mathcal{L}_{\frac{mq}{n},\infty}(H)$ are such that $$[A^{\frac{n}{2m}},B^{\frac{n}{m}}]\in(\mathcal{L}_{\frac{mq}{n},\infty})_0(H).$$
Applying Lemma \ref{CSZ main lemma} with $p=\frac{mq}{n}$ to the couple $(A^{\frac{n}{m}},B^{\frac{n}{m}}),$ we obtain
\begin{equation}\label{postdenis eq1}
B^qA^q-(A^{\frac{n}{2m}}B^{\frac{n}{m}}A^{\frac{n}{2m}})^{\frac{mq}{n}}\in(\mathcal{L}_{1,\infty})_0(H).
\end{equation}

Combining \eqref{postdenis eq0} and \eqref{postdenis eq1}, we write
$$(A^{\frac{ns}{2m}}B^{\frac{ns}{m}}A^{\frac{ns}{2m}})^{\frac{mq}{ns}}-
(A^{\frac{n}{2m}}B^{\frac{n}{m}}A^{\frac{n}{2m}})^{\frac{mq}{n}}
\in(\mathcal{L}_{1,\infty})_0(H).$$
Using Lemma \ref{bks lemma} with
$$X=(A^{\frac{ns}{2m}}B^{\frac{ns}{m}}A^{\frac{ns}{2m}})^{\frac{mq}{ns}},\quad Y=(A^{\frac{n}{2m}}B^{\frac{n}{m}}A^{\frac{n}{2m}})^{\frac{mq}{n}},\quad r=1,\quad t=\frac{ns}{mq},$$
we write
\begin{equation}\label{postdenis eq2}
A^{\frac{ns}{2m}}B^{\frac{ns}{m}}A^{\frac{ns}{2m}}-
(A^{\frac{n}{2m}}B^{\frac{n}{m}}A^{\frac{n}{2m}})^s\in(\mathcal{L}_{\frac{mq}{ns},\infty})_0(H).
\end{equation}

By assumption,
$$A^{\frac{n}{2m}}B^{\frac{n}{m}}A^{\frac{n}{2m}}-(A^{\frac1{2m}}B^{\frac1{m}}A^{\frac1{2m}})^n\in (\mathcal{L}_{\frac{mq}{n},\infty})_0(H).$$
Using Lemma \ref{bks lemma} with
$$X=A^{\frac{n}{2m}}B^{\frac{n}{m}}A^{\frac{n}{2m}},\quad Y=(A^{\frac1{2m}}B^{\frac1{m}}A^{\frac1{2m}})^n,\quad r=\frac{mq}{n},\quad t=s,$$
we write
\begin{equation}\label{postdenis eq3}
(A^{\frac{n}{2m}}B^{\frac{n}{m}}A^{\frac{n}{2m}})^s-(A^{\frac1{2m}}B^{\frac1{m}}A^{\frac1{2m}})^{ns}
\in(\mathcal{L}_{\frac{mq}{ns},\infty})_0(H).
\end{equation}
By assumption,
\begin{equation}\label{postdenis eq4}
A^{\frac{ns}{2m}}B^{\frac{ns}{m}}A^{\frac{ns}{2m}}-(A^{\frac{s}{2m}}B^{\frac{s}{m}}A^{\frac{s}{2m}})^n
\in (\mathcal{L}_{\frac{mq}{ns},\infty})_0(H).
\end{equation}
Combining \eqref{postdenis eq2}, \eqref{postdenis eq3} and \eqref{postdenis eq4}, we write
$$(A^{\frac{s}{2m}}B^{\frac{s}{m}}A^{\frac{s}{2m}})^n-(A^{\frac1{2m}}B^{\frac1{m}}A^{\frac1{2m}})^{ns}
\in(\mathcal{L}_{\frac{mq}{ns},\infty})_0(H).$$
Using Lemma \ref{bks lemma} with
$$X=(A^{\frac{s}{2m}}B^{\frac{s}{m}}A^{\frac{s}{2m}})^n,\quad Y=(A^{\frac1{2m}}B^{\frac1{m}}A^{\frac1{2m}})^{ns},\quad r=\frac{mq}{ns},\quad s=\frac1n,$$
we write
\begin{equation}\label{postdenis eq5}
A^{\frac{s}{2m}}B^{\frac{s}{m}}A^{\frac{s}{2m}}-(A^{\frac1{2m}}B^{\frac1{m}}A^{\frac1{2m}})^s
\in(\mathcal{L}_{\frac{mq}{s},\infty})_0(H).
\end{equation}

By assumption,
$$(A^{\frac1{2m}}B^{\frac1{m}}A^{\frac1{2m}})^m-A^{\frac12}BA^{\frac12}\in(\mathcal{L}_{q,\infty})_0(H).$$
Using Lemma \ref{bks lemma} with
$$X=(A^{\frac1{2m}}B^{\frac1{m}}A^{\frac1{2m}})^m,\quad Y=A^{\frac12}BA^{\frac12},\quad r=q,\quad t=\frac{s}{m},$$
we write
\begin{equation}\label{postdenis eq6}
(A^{\frac1{2m}}B^{\frac1{m}}A^{\frac1{2m}})^s-(A^{\frac12}BA^{\frac12})^{\frac{s}{m}}
\in(\mathcal{L}_{\frac{mq}{s},\infty})_0(H).
\end{equation}
Combining \eqref{postdenis eq5} and \eqref{postdenis eq6}, we write
$$A^{\frac{s}{2m}}B^{\frac{s}{m}}A^{\frac{s}{2m}}-(A^{\frac12}BA^{\frac12})^{\frac{s}{m}}
\in(\mathcal{L}_{\frac{mq}{s},\infty})_0(H).$$
By Lemma \ref{bks lemma},
$$(A^{\frac{s}{2m}}B^{\frac{s}{m}}A^{\frac{s}{2m}})^m-(A^{\frac12}BA^{\frac12})^s
\in(\mathcal{L}_{\frac{q}{s},\infty})_0(H).$$
By assumption,
$$(A^{\frac{s}{2m}}B^{\frac{s}{m}}A^{\frac{s}{2m}})^m-A^{\frac{s}{2}}B^sA^{\frac{s}{2}}
\in(\mathcal{L}_{\frac{q}{s},\infty})_0(H).$$
Combining the last two inclusions, we complete the proof.
\end{proof}

\subsection{Application of abstract estimate}\label{Application of abstract estimate}
We prove Theorem \ref{main example thm} in this subsection. To do that, we introduce the notion for differentiability of operators inspiring by Connes--Moscovici \cite[page 209]{CM1995}. See the following definition on $\mathbb{R}^{d}$. 
\begin{definition}\label{differentiable euclidean} 
Element $S\in\Pi_{\mathbb{R}^d}$ is called $m$ times differentiable if, for every second order differential operator $P$ with smooth compactly supported coefficients, $k$-th commutator
$$\underbrace{[P,[P,\cdots[P,S]\cdots]}_{k \text{ times}}$$
extends to a bounded operator from $W^{k,2}(\mathbb{R}^d)\to L_2(\mathbb{R}^d)$ for every $1\leq k\leq m.$
\end{definition}

The similar notion on compact manifold can be proposed in the following way.
\begin{definition}\label{differentiable manifold} 
Let $(X,G)$ be a compact manifold. Element $S\in\Pi_X$ is called $m$ times differentiable if, for every second order differential operator $P$ with smooth coefficients, $k$-th commutator
$$\underbrace{[P,[P,\cdots[P,S]\cdots]}_{k \text{ times}}$$
extends to a bounded operator from $W^{k,2}(X,{\rm vol}_{G})\to L_2(X,{\rm vol}_{G})$ for every $1\leq k\leq m.$
\end{definition}

Let us turn to the density of differentiable elements.  Notably, the algebra $\{g(\frac{\nabla}{\sqrt{\Delta}}):g\in C(\mathbb{S}^{d-1})\}$ can also be generated by the Riesz transforms $R_{j}=D_{j}\Delta^{-\frac{1}{2}}$, $j=1,\cdots,d$. We have below results.

\begin{lemma}\label{euclidean density lemma} 
Twice differentiable elements are dense in $\Pi_{\mathbb{R}^d}$. 
\end{lemma}
\begin{proof} 
Let $P$ be a second order differential operator with smooth compactly supported coefficients. Let $f\in C_{c}^{\infty}(\mathbb{R}^d)$. By \eqref{pseudo-differential op}, 
$$[P,M_f](1+\Delta)^{-\frac{1}{2}},\quad [P,[P,M_f]](1+\Delta)^{-1}\in\Psi^{0}(\mathbb{R}^{d}).$$
Thus, the commutators on the left hand side are bounded, and therefore $M_f$ is twice differentiable.

Let's next show that Riesz transforms $\{R_j\}_{j=1}^d$ are twice differentiable. It suffices to show
$$[P,R_{j}](1+\Delta)^{-\frac{1}{2}},\quad [P,[P,R_{j}]](1+\Delta)^{-1}\in B(L_2(\mathbb{R}^d)).$$
The boundedness for the first operator can be obtained by adapting the proof of Lemma \ref{bd in second sum}. For the second one, write 
$$[P,[P,R_{j}]]=[P,[P,D_{j}(1+\Delta)^{-\frac{1}{2}}]]
+[P,[P,D_{j}(\Delta^{-\frac{1}{2}}-(1+\Delta)^{-\frac{1}{2}})]].$$
Clearly, $P\in\Psi^{2}(\mathbb{R}^{d})$ as $P$ has smooth compactly supported coefficients. By \eqref{pseudo-differential op},
$$[P,[P,D_{j}(1+\Delta)^{-\frac{1}{2}}]](1+\Delta)^{-1}\in\Psi^{0}(\mathbb{R}^{d}).$$
The boundedness of this operator follows from \cite[Theorem 1,page 235]{Stein-Murphy1993}. Observe that
$$A:=D_{j}(\Delta^{-\frac{1}{2}}-(1+\Delta)^{-\frac{1}{2}})
=(1+\Delta)^{-1}\cdot\frac{R_{j}\sqrt{1+\Delta}}{\sqrt{1+\Delta}+\sqrt{\Delta}},$$
and 
$$[P,[P,A]]=P^{2}A-2PAP+AP^{2}.$$
Write
$$P^{2}A(1+\Delta)^{-1}=P^{2}(1+\Delta)^{-2}
\cdot\frac{R_{j}\sqrt{1+\Delta}}{\sqrt{1+\Delta}+\sqrt{\Delta}}.$$
On the right hand side of the equality, the first factor belongs to $\Psi^{0}(\mathbb{R}^{d})$ due to $P^{2}\in\Psi^{4}(\mathbb{R}^{d})$ and \eqref{pseudo-differential op}, and therefore is bounded on $L_{2}(\mathbb{R}^{d})$. The boundedness for the second one is obvious. Similarly, the operators
$$PAP(1+\Delta)^{-1}=P(1+\Delta)^{-1}
\cdot\frac{R_{j}\sqrt{1+\Delta}}{\sqrt{1+\Delta}+\sqrt{\Delta}}\cdot P(1+\Delta)^{-1}$$
and
$$AP^{2}(1+\Delta)^{-1}=\frac{R_{j}\sqrt{1+\Delta}}{\sqrt{1+\Delta}+\sqrt{\Delta}}
\cdot(1+\Delta)^{-1}P^{2}(1+\Delta)^{-1}$$
are bounded on $L_{2}(\mathbb{R}^{d})$. In a sum, $[P,[P,A]](1+\Delta)^{-1}$ is bounded on $L_{2}(\mathbb{R}^{d})$.
Moreover, a polynomial $S$ generated by $\{M_f: f\in C_{c}^{\infty}(\mathbb{R}^d)\}$ and $\{R_j\}_{j=1}^d$ is twice differentiable.

For $S\in\Pi_{\mathbb{R}^{d}}$, pick a sequence $\{S_n\}_{n\geq0}$ such that each $S_n$ belongs to $\ast$-algebra generated by $\{M_f: f\in C_0(\mathbb{R}^d)\}$ and by $\{R_j\}_{j=1}^d$, and such that $S_n\to S$ in the uniform norm.  Since $C_{c}^{\infty}(\mathbb{R}^{d})$ is dense in $C_{0}(\mathbb{R}^{d})$ in the uniform norm, there is a sub-sequence $\{S_{n_{k}}\}_{k\geq0}$ generated by $\{M_f: f\in C_{c}^{\infty}(\mathbb{R}^d)\}$ and $\{R_j\}_{j=1}^d$, and $S_{n_{k}}\to S$ in the uniform norm. By the previous arguments, $S_{n_{k}}$ is twice differentiable. The assertion follows.
\end{proof}

\begin{lemma}\label{local density lemma} 
Let $(X,G)$ be a compact Riemannian manifold. Twice differentiable elements compactly supported in the chart $(\U_i,h_i)$ are dense in the set of elements of $\Pi_X$ compactly supported in the chart $(\U_i,h_i).$
\end{lemma}
\begin{proof} Let $(\mathcal{U}_i,h_i)$ be a chart and let $S\in\Pi_X$ be compactly supported  in $\mathcal{U}_i.$ The operator $W_{h_i}SW_{h_i}^{\ast}\in\Pi_{\mathbb{R}^d}$ is compactly supported in $h_i(\mathcal{U}_i).$ Fix $\phi_i\in C^{\infty}_c(h_i(\mathcal{U}_i))$ such that
$$W_{h_i}SW_{h_i}^{\ast}=W_{h_i}SW_{h_i}^{\ast}\cdot M_{\phi_i}=M_{\phi_i}\cdot W_{h_i}SW_{h_i}^{\ast}.$$

By Lemma \ref{euclidean density lemma}, there exists a sequence $\{S_n\}_{n\geq1}\subset\Pi_{\mathbb{R}^d}$ such that (i) $S_n$ is twice differentiable for every $n\geq1$ and (ii) $S_n\to W_{h_i}SW_{h_i}^{\ast}$ in the uniform norm. Since $M_{\phi_i}$ is\footnote{This is established in the proof of Lemma \ref{euclidean density lemma}.} twice differentiable, it follows that so is $M_{\phi_i}S_nM_{\phi_i}$ for every $n\geq1.$ We have
$$M_{\phi_i}S_nM_{\phi_i}\to M_{\phi_i}\cdot W_{h_i}SW_{h_i}^{\ast}\cdot M_{\phi_i}=W_{h_i}SW_{h_i}^{\ast}$$
in the uniform norm. Replacing $S_n$ with $M_{\phi_i}S_nM_{\phi_i}$, we may assume without loss of generality that each $S_n$ is compactly supported in $h_i(\mathcal{U}_i).$ In this case, $W_{h_i}^{\ast}S_nW_{h_i}\to S$ in the uniform norm and therefore it suffices to show that $W_{h_i}^{\ast}S_nW_{h_i}\in\Pi_X$ is twice differentiable for every $n\geq 1.$

Let $P$ be a second order differential operator with smooth coefficients on $X$. Let us prove that 
\begin{align}\label{WSW condition1}
[W_{h_{i}}^{*}S_{n}W_{h_{i}},P]:W^{1,2}(X,{\rm vol}_{G})\rightarrow L_{2}(X,{\rm vol}_{G})
\end{align}
and 
\begin{align}\label{WSW condition2}
[[W_{h_{i}}^{*}S_{n}W_{h_{i}},P],P]:W^{2,2}(X,{\rm vol}_{G})\rightarrow L_{2}(X,{\rm vol}_{G})
\end{align}
are bounded. We begin by considering \eqref{WSW condition1}.
Denote $\psi_{i}=\phi_{i}\circ h_{i}$ and write
\begin{equation}\label{split term}
\begin{split}
&[W_{h_{i}}^{*}S_{n}W_{h_{i}},P]\\
&=W_{h_{i}}^{*}S_{n}W_{h_{i}}[M_{\psi_{i}},P]+[W_{h_{i}}^{*}S_{n}W_{h_{i}},M_{\psi_{i}}PM_{\psi_{i}}]
+[P,M_{\psi_{i}}]W_{h_{i}}^{*}S_{n}W_{h_{i}}.
\end{split}
\end{equation}
The first and third factors on the right hand side of equality are bounded from $W^{1,2}(X,{\rm vol}_{G})$ to $L_{2}(X,{\rm vol}_{G})$ since $W_{h_{i}}^{*}S_{n}W_{h_{i}}$ is bounded on $L_{2}(X,{\rm vol}_{G})$ and $[M_{\psi_{i}},P]:W^{1,2}(X,{\rm vol}_{G})\to L_{2}(X,{\rm vol}_{G})$ is bounded. 
For the second factor, write 
$$[W_{h_{i}}^{*}S_{n}W_{h_{i}},M_{\psi_{i}}PM_{\psi_{i}}]=W_{h_{i}}^{*}[S_{n},P_{i}]W_{h_{i}}$$
where $P_{i}=W_{h_{i}}M_{\psi_{i}}PM_{\psi_{i}}W_{h_{i}}^{*}$ is a second order differential operator with smooth coefficients compactly supported in $\Omega_{i}$. 
Pick $\eta_{i}\in C_{c}^{\infty}(\Omega_{i})$ such that 
$$[S_{n},P_{i}]=M_{\eta_{i}}[S_{n},P_{i}]M_{\eta_{i}}.$$
Since $S_{n}$ is twice differentiable, it follows that $[S_{n},P_{i}]:W^{1,2}(\Omega_{i})\rightarrow L_{2}(\Omega_{i})$ is bounded.
Clearly, 
$$M_{\eta_{i}}:W^{1,2}(\Omega_{i},{\rm vol}_{g_{i}})\rightarrow W^{1,2}(\Omega_{i})$$
and 
$$M_{\eta_{i}}:L_{2}(\Omega_{i})\rightarrow L_{2}(\Omega_{i},{\rm vol}_{g_{i}})$$ 
are bounded. 
In a sum, $[S_{n},P_{i}]:W^{1,2}(\Omega_{i},{\rm vol}_{g_{i}})\rightarrow L_{2}(\Omega_{i},{\rm vol}_{g_{i}})$ is bounded and therefore the assertion for \eqref{WSW condition1} holds.

Turn to investigate \eqref{WSW condition2}. Firstly,
$$[W_{h_{i}}^{*}S_{n}W_{h_{i}}[M_{\psi_{i}},P],P]=W_{h_{i}}^{*}S_{n}W_{h_{i}}[[M_{\psi_{i}},P],P]
+[W_{h_{i}}^{*}S_{n}W_{h_{i}},P][M_{\psi_{i}},P].$$
The first factor on the right hand side is bounded from $W^{2,2}(X,{\rm vol}_{G})$ to $L_{2}(X,{\rm vol}_{G})$ since $[[M_{\psi_{i}},P],P]:W^{2,2}(X,{\rm vol}_{G})\to L_{2}(X,{\rm vol}_{G})$ is bounded and $W_{h_{i}}^{*}S_{n}W_{h_{i}}$ is bounded on $L_{2}(X,{\rm vol}_{G})$. 
The second one is bounded from $W^{2,2}(X,{\rm vol}_{G})$ to $L_{2}(X,{\rm vol}_{G})$ since operator $[M_{\psi_{i}},P]:W^{2,2}(X,{\rm vol}_{G})\to W^{1,2}(X,{\rm vol}_{G})$ is bounded and $[W_{h_{i}}^{*}S_{n}W_{h_{i}},P]:W^{1,2}(X,{\rm vol}_{G})\to L_{2}(X,{\rm vol}_{G})$ is bounded by \eqref{WSW condition1}. Similarly,
$$[[P,M_{\psi_{i}}]W_{h_{i}}^{*}S_{n}W_{h_{i}},P]=+[P,M_{\psi_{i}}][W_{h_{i}}^{*}S_{n}W_{h_{i}},P]
+[[P,M_{\psi_{i}}],P]W_{h_{i}}^{*}S_{n}W_{h_{i}}$$
is bounded from $W^{2,2}(X,{\rm vol}_{G})$ to $L_{2}(X,{\rm vol}_{G})$.
Let $A_{i}=[W_{h_{i}}^{*}S_{n}W_{h_{i}},M_{\psi_{i}}PM_{\psi_{i}}]$. We have $A_{i}=M_{\psi_{i}}A_{i}M_{\psi_{i}}$. Moreover,
$$[A_{i},P]=A_{i}[M_{\psi_{i}},P]+[A_{i},M_{\psi_{i}}PM_{\psi_{i}}]
+[P,M_{\psi_{i}}]A_{i}.$$
The first factor is bounded from $W^{2,2}(X,{\rm vol}_{G})$ to $L_{2}(X,{\rm vol}_{G})$ since $[M_{\psi_{i}},P]$ is bounded from $W^{2,2}(X,{\rm vol}_{G})$ to $W^{1,2}(X,{\rm vol}_{G})$ and $A_{i}$ is bounded from $W^{1,2}(X,{\rm vol}_{G})$ to $L_{2}(X,{\rm vol}_{G})$. The third factor is similar. For the second one, write
$$[A_{i},M_{\psi_{i}}PM_{\psi_{i}}]=W_{h_{i}}^{*}[[S_{n},P_{i}],P_{i}]W_{h_{i}}$$
where $P_{i}=W_{h_{i}}M_{\psi_{i}}PM_{\psi_{i}}W_{h_{i}}^{*}$. Since $S_{n}$ is twice differentiable, it follows that $[[S_{n},P_{i}],P_{i}]:W^{2,2}(\Omega_{i},{\rm vol}_{g_{i}})\rightarrow L_{2}(\Omega_{i},{\rm vol}_{g_{i}})$ is bounded by adapting the arguments in the last second paragraph. Lastly, by \eqref{split term}, the assertion for \eqref{WSW condition2} holds.
This completes the proof.
\end{proof}

\begin{lemma}\label{main example density lemma} 
Let $(X,G)$ be a compact Riemannian manifold. Twice differentiable elements are dense in $\Pi_X$ in the uniform norm.
\end{lemma}
\begin{proof}
Let $\{\phi_{n}\}_{n=1}^{N}\subset C_{c}^{\infty}(X)$ be a partition of unity subordinate to the atlas of $X.$ For $S\in\Pi_X,$ we write 
$$S=\sum_{n=1}^{N}M_{\phi_{n}^{\frac{1}{2}}}SM_{\phi_{n}^{\frac{1}{2}}}+
\sum_{n=1}^{N}[S,M_{\phi_{n}^{\frac{1}{2}}}]M_{\phi_{n}^{\frac{1}{2}}}.$$

The approximation for summands in the first sum on the right hand side are given by Lemma \ref{local density lemma}. Observe that the summands in the second sum on the right hand side are compact operators. Since integral operators with smooth integral kernels are twice differentiable, compact and dense among compact operators in the uniform norm, the assertion follows.
\end{proof}

\begin{lemma}\label{differentiable-bd} 
Let $(X,G)$ be a compact Riemannian manifold. If $S\in\Pi_X$ is twice differentiable, then $[S,(1+\Delta_G)^{\frac12}]$ is bounded.
\end{lemma}
\begin{proof} 
Since $X$ is compact, the assumption can be equivalently restated as the boundedness of operators
$$[\Delta_G,S](1+\Delta_G)^{-\frac12},[\Delta_G,[\Delta_G,S]](1+\Delta_G)^{-1}$$
on $L_2(X,{\rm vol}_G).$ Recall that 
$$(1+\Delta_G)^{-\frac{1}{2}}=
\frac{1}{\pi}\int_{0}^{\infty}(1+\lambda+\Delta_{G})^{-1}\lambda^{-\frac{1}{2}}d\lambda$$
where the convergence for the integral is understood in the operator norm. Write
\begin{align}\label{D^12_G decomposition}
[(1+\Delta_G)^{\frac{1}{2}},S]=[(1+\Delta_G)(1+\Delta_G)^{-\frac{1}{2}},S]
=\frac{1}{\pi}\int_{0}^{\infty}[\frac{1+\Delta_{G}}{1+\lambda+\Delta_{G}},S]
\lambda^{-\frac{1}{2}}d\lambda.
\end{align}
However,
\begin{align*}
[\frac{1+\Delta_{G}}{1+\lambda+\Delta_{G}},S]
&=[\Delta_{G},S](1+\Delta_{G})^{-\frac{1}{2}}\cdot
\Big(\frac{(1+\Delta_{G})^{\frac{1}{2}}}{1+\lambda+\Delta_{G}}-
\frac{(1+\Delta_{G})^{\frac{3}{2}}}{(1+\lambda+\Delta_{G})^{2}}\Big)+\\
&+\frac{\lambda}{1+\lambda+\Delta_{G}}\cdot[\Delta_{G},[\Delta_{G},S]](1+\Delta_{G})^{-1}
\cdot\frac{1+\Delta_{G}}{(1+\lambda+\Delta_{G})^{2}}\\
&=:X_{\lambda}+Y_{\lambda}.
\end{align*}
By functional calculus,
$$\int_{0}^{\infty}\Big(\frac{(1+\Delta_{G})^{\frac{1}{2}}}{1+\lambda+\Delta_{G}}-
\frac{(1+\Delta_{G})^{\frac{3}{2}}}{(1+\lambda+\Delta_{G})^{2}}\Big)\lambda^{-\frac{1}{2}}d\lambda
=\frac{\pi}{2}.$$
It follows that
$$\frac{1}{\pi}\int_{0}^{\infty}X_{\lambda}\cdot\lambda^{-\frac{1}{2}}d\lambda=
\frac{1}{2}[\Delta_{G},S](1+\Delta_{G})^{-\frac{1}{2}}$$
is bounded on $L_2(X,{\rm vol}_G)$ by the assumption.  
Additionally, by spectral theorem,
$$\|\frac{\lambda}{1+\lambda+\Delta_{G}}\|_{B(L_2(X,{\rm vol}_G))}=\frac{\lambda}{1+\lambda}$$
and
\begin{equation*}
\|\frac{1+\Delta_{G}}{(1+\lambda+\Delta_{G})^{2}}\|_{B(L_2(X,{\rm vol}_G))}
=\left\{
\begin{array}{cl}
\frac{1}{(1+\lambda)^{2}},\quad& 0<\lambda\leq1,\\
\frac{1}{4\lambda},\quad& \lambda\geq1.
 \end{array}\right.
\end{equation*}
Thus,
\begin{align*}
\|\int_{0}^{\infty}Y_{\lambda}\cdot\lambda^{-\frac{1}{2}}d\lambda\|_{B(L_2(X,{\rm vol}_G))}
\leq\frac{7}{6}\|[\Delta_{G},[\Delta_{G},S]](1+\Delta_{G})^{-1}\|_{B(L_2(X,{\rm vol}_G))}.
\end{align*}
Taking the arguments for $X_{\lambda},Y_{\lambda}$, equality \eqref{D^12_G decomposition} and the assumption for $S$ into a combination, the assertion follows.
\end{proof}

Next, we verify the assumptions of Theorem \ref{nearly commuting powers} on manifold $(X,G)$.
\begin{lemma}\label{main example verification lemma} Let $(X,G)$ be a compact Riemannian manifold and let $\Delta_G$ be the Laplace-Beltrami operator. If $0\leq S\in\Pi_X,$ then
$A=S$ and $B=(1+\Delta_G)^{-1}$
satisfy the assumptions in Theorem \ref{nearly commuting powers} with $q=\frac{d}{2}.$
\end{lemma}
\begin{proof} Let $\{S_n\}_{n\geq1}\subset \Pi_X$ be a sequence of twice differentiable elements such that $S_n\to S^{\epsilon}$ in the uniform norm.

Fix $m>\delta$ and write
\begin{align*}
[S_n,B^{\delta}]=[S_n,(B^{\frac{\delta}{2m}})^{2m}]
&=\sum_{k=0}^{2m-1}(B^{\frac{\delta}{2m}})^k\cdot [S_n,B^{\frac{\delta}{2m}}] \cdot (B^{\frac{\delta}{2m}})^{2m-1-k}\\
&=-\sum_{k=0}^{2m-1}(B^{\frac{\delta}{2m}})^{k+1}\cdot [S_n,B^{-\frac{\delta}{2m}}]\cdot (B^{\frac{\delta}{2m}})^{2m-k}.
\end{align*}
As $S_n$ is twice differentiable, the operator $[S_n,B^{-\frac12}]$ is bounded due to Lemma \ref{differentiable-bd}. By Birman-Koplienko-Solomyak inequality \cite[Theorem 1.3.4]{LMSZvol2}, $$\int_{0}^{t}\mu(s,[S_n,B^{-\frac{\delta}{2m}}])ds\leq c_{\delta,m}\|S_{n}\|_{B(L_2(X,{\rm vol}_G))}^{1-\frac{\delta}{m}}
\int_{0}^{t}\mu(s,|[S_n,B^{-\frac{1}{2}}]|^{\frac{\delta}{m}})ds,\quad t>0.$$
Thus, the commutator $[S_n,B^{-\frac{\delta}{2m}}]$ is bounded. It follows now from H\"older inequality that
$$[S_n,B^{\delta}]\in\mathcal{L}_{\frac{dm}{(2m+1)\delta},\infty}(L_2(X,{\rm vol}_G))\subset (\mathcal{L}_{\frac{d}{2\delta}})_0(L_2(X,{\rm vol}_G)).$$
Since $[S_n,B^{\delta}]\to [A^{\epsilon},B^{\delta}]\in \mathcal{L}_{\frac{d}{2\delta}}(L_2(X,{\rm vol}_G))$ as $n\to\infty,$ the assertion follows from Lemma \ref{main example density lemma}. The proof is complete.
\end{proof}

\begin{proof}[Proof of Theorem \ref{main example thm}] It is established in Lemma \ref{main example verification lemma} that the operators
$$A=S,\quad B=(1+\Delta_G)^{-1}$$
satisfy the assumptions in Theorem \ref{nearly commuting powers} with $q=\frac{d}{2}.$ Applying Theorem \ref{nearly commuting powers}, we complete the proof.
\end{proof}

\section{Asymptotics for operators localised in a chart}\label{sec-6}
This section is devoted to the proof of Theorem \ref{asymptotic local}, which establishes a local version of the result stated in Theorem \ref{asymptotic on manifold}.
\subsection{Good atlas on compact Riemannian manifold}\label{Good atlas}
To show Theorem \ref{asymptotic local}, we introduce the definition of good atlas on compact manifold.
\begin{definition}\label{good atlas def} Let $(X,G)$ be a compact Riemannian manifold. A finite atlas $\{(\U_i,h_i)\}_{i\in\mathbb{I}}$ is called good if, for every $i\in\mathbb{I},$ the function $g_i:h_i(\U_i)\to {\rm GL}_+(d,\mathbb{R})$ allows a smooth extension $g_i:\mathbb{R}^d\to{\rm GL}_+(d,\mathbb{R})$ which is constant outside some ball.
\end{definition}

Towards constructing good atlas on compact Riemannian manifolds, the following smoothness result is established as a start. 
\begin{lemma}\label{log composition lemma}
Let $\Omega$ be a domain in $\mathbb{R}^d$ and let $g:\Omega\to {\rm GL}_+(d,\mathbb{R})$ be a smooth function. The function $\log\circ g:\Omega\to M_d(\mathbb{R})$ is also smooth where $M_d(\mathbb{R})$ denotes the space of all $d\times d$ real matrices.
\end{lemma}
\begin{proof} 
Fix $t_0\in\Omega.$ Since $g(t_0)$ is a positive matrix, its spectrum is positive. Let $\lambda_{\text{min}}(t_0)$ and $\lambda_{\text{max}}(t_0)$ denote its smallest and largest eigenvalues, respectively. We have
\begin{align*}
\|g(t_0)\|_{M_d(\mathbb{R})}=\lambda_{\text{max}}(t_0), \quad 
\|g(t_0)^{-1}\|_{M_d(\mathbb{R})}=\frac{1}{\lambda_{\text{min}}(t_0)},
\end{align*}
where $\|\cdot\|_{M_d(\mathbb{R})}$ denotes the operator norm on $M_d(\mathbb{R})$. 
Therefore,
$$\|g(t_0)^{-1}\|_{M_d(\mathbb{R})}^{-1}\leq g(t_0)\leq\|g(t_0)\|_{M_d(\mathbb{R})}.$$
Choose $\epsilon>0$ such that $$\|g(t)-g(t_0)\|_{M_d(\mathbb{R})}<\frac12\|g(t_0)^{-1}\|_{M_d(\mathbb{R})}^{-1},\quad |t-t_0|<\epsilon.$$ 
We have
$$\frac12\|g(t_0)^{-1}\|_{M_d(\mathbb{R})}^{-1}\leq g(t)\leq \frac12\|g(t_0)^{-1}\|_{M_d(\mathbb{R})}^{-1}+\|g(t_0)\|_{M_d(\mathbb{R})},\quad |t-t_0|<\epsilon.$$
Choose $\omega\in C^{\infty}(\mathbb{R})$ such that 
$$\omega(u)=\log(u),\quad u\in[\frac12\|g(t_0)^{-1}\|_{M_d(\mathbb{R})}^{-1},\frac12\|g(t_0)^{-1}\|_{M_d(\mathbb{R})}^{-1}+\|g(t_0)\|_{M_d(\mathbb{R})}].$$
The mapping $\omega\circ g$ is obviously smooth on $\Omega$. Note that $\log\circ g=\omega\circ g$ on $B(t_0,\epsilon).$ Hence, $\log\circ g$ is smooth on $B(t_0,\epsilon).$ Since $t_0\in\Omega$ is arbitrary, the assertion follows.
\end{proof}

A compact manifold admits a finite atlas. This allows us to derive a more refined atlas in the setting of compact Riemannian manifolds.
\begin{proposition}\label{good atlas prop}
Every compact Riemannian manifold allows a good atlas.
\end{proposition}
\begin{proof}
Let $(X, G)$ be a compact Riemannian manifold with atlas $\{(\U_i,h_i)\}_{i\in\mathbb{I}}.$ For every open set $D$ compactly supported in $h_i(\U_i)$ (i.e. the closure $\overline{D}$ is compact and $\overline{D} \subset h_i(\mathcal{U}_i)$), choose a function $\phi_D\in C^{\infty}_c(h_i(\U_i))$ such that $\phi_D=1$ on $D.$ Since $X$ is a Riemannian manifold, the function $g_i:h_{i}(\mathcal{U}_{i})\to{\rm GL}_+(d,\mathbb{R})$ is well-defined and smooth. By Lemma \ref{log composition lemma}, the function $\log\circ g_i:h_i(\U_i)\to M_d(\mathbb{R})$ is smooth on $h_i(\U_i).$ Hence, $g_i:D\to{\rm GL}_+(d,\mathbb{R})$ allows an extension $g_{i,D}:\mathbb{R}^d\to {\rm GL}_+(d,\mathbb{R})$ given by the formula
$$g_{i,D}:x\to \exp(\phi_{D}(x)\cdot \log(g_i(x))),\quad x\in\mathbb{R}^d.$$
Clearly, $g_{i,D}$ is a smooth extension of $g_i|_{D}$ and it is constant outside some ball.

Consider now the new atlas \footnote{Here, $D\Subset h_i(\U_i)$ means "$D$ is compactly supported in $h_i(\U_i)$".} $\{\{(h_i^{-1}(D),h_i)\}_{D\Subset h_i(\U_i)}\}_{i\in\mathbb{I}}.$ In this atlas, Riemannian metric $V$ is given by the formula $V_{(i,D)}=G_i|_{h_i^{-1}(D)}$ and the coordinate mapping $h_{(i,D)}$ is $h_{i}$ related to  $h_i^{-1}(D)$. For every index $(i,D),$ we have $V_{(i,D)}\circ h_{(i,D)}^{-1}=g_{i,D}|_D=g_i|_{D}.$ Finally, by compactness on $X$, there exists a finite sub-atlas and it is good. The proof is complete.
\end{proof}

\subsection{Local case for $p=\frac{d}{2}$}
This subsection aims to show Theorem \ref{asymptotic local}. We need the necessary lemmas from below.
\begin{lemma}\label{boundedness on manifold}
Let $(X,G)$ be a compact Riemannian manifold and let $\Delta_G$ be the Laplace-Beltrami operator. Let $(\U_i,h_i)$ be a chart of $X$. For $\phi\in C^{\infty}_c(\U_i),$ the operator $[\Delta_G,M_{\phi}](1+\Delta_G)^{-\frac12}$ is bounded. 
\end{lemma}
\begin{proof}
Clearly, $M_{\phi}\in\Psi^{0}(X)$, $\Delta_{G}\in\Psi^{2}(X)$ and $(1+\Delta_G)^{-\frac12}\in\Psi^{-1}(X)$. We have $[\Delta_{G},M_{\phi}]\in\Psi^{1}(X)$ by \cite[Remark 5.3.3]{RT2009-book}. The composition of pseudo-differential operators gives that $[\Delta_{G},M_{\phi}](1+\Delta_{G})^{-\frac{1}{2}}\in\Psi^{0}(X).$ 
Since $X$ is a compact manifold, the assertion follows from \cite[Theorem 5.2.23]{RT2009-book}. 
\end{proof}

\begin{lemma}\label{yuri lemma} Let $(X,G)$ be a compact Riemannian manifold and let $\Delta_G$ be the Laplace-Beltrami operator. Let $(\U_i,h_i)$ be a chart in a good (in the sense of Definition \ref{good atlas def}) atlas. For every $\phi\in C^{\infty}_c(\U_i),$ we have
$$M_{\phi}(1+\Delta_G)^{-1}M_{\phi}-W_{h_i}^{\ast}M_{\phi\circ h_i^{-1}}(1+\Delta_{g_i})^{-1}M_{\phi\circ h_i^{-1}}W_{h_i}\in\mathcal{L}_{\frac{d}{3},\infty}(L_2(X,{\rm vol}_G)).$$
\end{lemma}
\begin{proof} We prove an equivalent assertion
$$W_{h_i}M_{\phi}(1+\Delta_G)^{-1}M_{\phi}W_{h_i^{\ast}}-M_{\phi\circ h_i^{-1}}(1+\Delta_{g_i})^{-1}M_{\phi\circ h_i^{-1}}\in\mathcal{L}_{\frac{d}{3},\infty}(L_2(\mathbb{R}^d,{\rm vol}_{g_i})).$$
Denote the operator on the left hand side in the latter inclusion by $LHS.$

Note that $LHS$ sends $L_2(\mathbb{R}^d,{\rm vol}_{g_i})$ into $W^{2,2}(\mathbb{R}^d,{\rm vol}_{g_i}),$ which is the domain of $\Delta_{g_i}.$ We have
\begin{align*}
&(1+\Delta_{g_i})\cdot LHS=\\
&=W_{h_i}(1+\Delta_G)M_{\phi}(1+\Delta_G)^{-1}M_{\phi}W_{h_i^{\ast}}-(1+\Delta_{g_i})M_{\phi\circ h_i^{-1}}(1+\Delta_{g_i})^{-1}M_{\phi\circ h_i^{-1}}\\
&=W_{h_i}[\Delta_G,M_{\phi}](1+\Delta_G)^{-1}M_{\phi}W_{h_i^{\ast}}-[\Delta_{g_i},M_{\phi\circ h_i^{-1}}](1+\Delta_{g_i})^{-1}M_{\phi\circ h_i^{-1}}.
\end{align*}
In other words,
$$LHS=(1+\Delta_{g_i})^{-1}\cdot\Big(W_{h_i}ABW_{h_i}^{\ast}-CD\Big),$$
where
$$A=[\Delta_G,M_{\phi}](1+\Delta_G)^{-\frac12},\quad B=(1+\Delta_G)^{-\frac12}M_{\phi},$$
$$C=[\Delta_{g_i},M_{\phi\circ h_i^{-1}}](1+\Delta_{g_i})^{-\frac12},\quad D=(1+\Delta_{g_i})^{-\frac12}M_{\phi\circ h_i^{-1}}.$$

Choose $\psi\in C^{\infty}_c(\mathbb{R}^d)$ such that $\psi=1$ on and near ${\rm supp}(\phi\circ h_i^{-1}).$ It follows that
\begin{equation}\label{yl eq0}
LHS=(1+\Delta_{g_i})^{-1}M_{\psi}\cdot\Big(W_{h_i}ABW_{h_i}^{\ast}-CD\Big).
\end{equation}

By Lemma \ref{boundedness on manifold} and Lemma \ref{special cwikel manifold},
$$A\in\mathcal{L}_{\infty}(L_2(X,{\rm vol}_G)),\quad B\in \mathcal{L}_{d,\infty}(L_2(X,{\rm vol}_G)).$$
Hence,
\begin{equation}\label{yl eq1}
W_{h_i}ABW_{h_i}^{\ast}\in\mathcal{L}_{d,\infty}(L_2(X,{\rm vol}_G)).
\end{equation}

On the other hand, consider $M_{{\rm det}^{\frac14}(g_{i})}:L_2(\mathbb{R}^d,{\rm vol}_{g_{i}})\to L_2(\mathbb{R}^d)$ the unitary mapping. We have
$$M_{{\rm det}^{-\frac14}(g_{i})}[\Delta_{g_i},M_{\phi\circ h_i^{-1}}](1+\Delta_{g_i})^{-\frac12}M_{{\rm det}^{\frac14}(g_{i})}
=[P_{i},M_{\phi\circ h_i^{-1}}](1+P_{i})^{-\frac12}$$
where $P_{i}=M_{{\rm det}^{-\frac14}(g_{i})}\Delta_{g_i}M_{{\rm det}^{\frac14}(g_{i})}$ is a positive, self-adjoint, elliptic differential operator whose coefficients are smooth and also are constant outside some ball on $\mathbb{R}^{d}$. By Lemma \ref{P commutator bd-integer}, 
$$[P_{i},M_{\phi\circ h_i^{-1}}](1+P_{i})^{-\frac12}\in\mathcal{L}_{\infty}(L_2(\mathbb{R}^d))$$ 
and therefore $$C\in\mathcal{L}_{\infty}(L_2(\mathbb{R}^d,{\rm vol}_{g_i})).$$
Moreover, write $$M_{{\rm det}^{-\frac14}(g_{i})}(1+\Delta_{g_i})^{-\frac12}M_{\phi\circ h_i^{-1}}M_{{\rm det}^{\frac14}(g_{i})}=(1+P_{i})^{-\frac12}M_{\phi\circ h_i^{-1}}.$$
The Lemma \ref{P-Mf compactness} gives that
$$(1+P_{i})^{-\frac12}M_{\phi\circ h_i^{-1}}\in\mathcal{L}_{d,\infty}(L_2(\mathbb{R}^d))$$
and therefore
$$D\in\mathcal{L}_{d,\infty}(L_2(\mathbb{R}^d,{\rm vol}_{g_i})).$$
Hence,
\begin{equation}\label{yl eq2}
CD\in\mathcal{L}_{d,\infty}(L_2(\mathbb{R}^d,{\rm vol}_{g_i})).
\end{equation}

Note that
\begin{equation}\label{yl eq3} (1+\Delta_{g_i})^{-1}M_{\psi}\in\mathcal{L}_{d,\infty}(L_2(\mathbb{R}^d,{\rm vol}_{g_i})).
\end{equation}

Combining \eqref{yl eq0}, \eqref{yl eq1}, \eqref{yl eq2} and \eqref{yl eq3}, we complete the proof.
\end{proof}

\begin{lemma}\label{manifold compactness of commutators} For every $\phi\in C^{\infty}(X),$ we have
$$[M_{\phi},(1+\Delta_{G})^{-1}]\in \mathcal{L}_{\frac{d}{3},\infty}(L_2(X,{\rm vol}_G)).$$
\end{lemma}
\begin{proof} 
We have
\begin{align*}
[M_{\phi},(1+\Delta_{G})^{-1}]
&=(1+\Delta_{G})^{-1}[\Delta_G,M_{\phi}](1+\Delta_{G})^{-1}\\
&=(1+\Delta_{G})^{-1}\cdot [\Delta_G,M_{\phi}](1+\Delta_{G})^{-\frac12}\cdot(1+\Delta_{G})^{-\frac12}.
\end{align*}
The first factor belongs to $\mathcal{L}_{\frac{d}{2},\infty}(L_2(X,{\rm vol}_G)),$ the second factor is bounded due to Lemma \ref{boundedness on manifold}, the last factor belongs to $\mathcal{L}_{d,\infty}(L_2(X,{\rm vol}_G)).$ The assertion follows from H\"older inequality.
\end{proof}

The following theorem is a local version of the Theorem \ref{asymptotic on manifold} with $p=\frac{d}{2}$.
\begin{theorem}\label{asymptotic local} 
Let $(X,G)$ be a compact Riemannian manifold and let $\Delta_G$ be the Laplace-Beltrami operator. Let $(\U_i,h_i)$ be a chart in a good (in the sense of Definition \ref{good atlas def}) atlas. If $S\in\Pi_{X}$ is compactly supported in $\U_i,$ then
$$\lim_{t\rightarrow\infty}t^{\frac{2}{d}}\mu(t,S(1+\Delta_{G})^{-1})
=c_d\Big\|\sym_{X}(S)\Big\|_{L_{\frac{d}{2}}(T^{\ast}X,e^{-q_{G}}d\lambda)}$$
where $c_{d}$ is given in Theorem \ref{asymptotic on Rd}.
\end{theorem}

Before proving Theorem \ref{asymptotic local}, recall a fact that if $H_1$ and $H_2$ are Hilbert spaces and if $W:H_1\to H_2$ is a partial isometry, then
$$\mu_{B(H_2)}(WTW^{\ast})=\mu_{B(H_1)}(T),\quad T\in B(H),\quad T=T(W^{\ast}W)=(W^{\ast}W)T.$$	
Let us continue to show Theorem \ref{asymptotic local}.
\begin{proof}[Proof of Theorem \ref{asymptotic local}]
Fix $\phi\in C^{\infty}_c(\U_i)$ such that $S=M_{\phi}S=SM_{\phi}.$ Applying the above fact to the operator
$$T=S(1+\Delta_{G})^{-1}M_{\phi},$$
we conclude
\begin{equation}\label{alp eq0}
\mu_{B(L_2(\mathbb{R}^d,{\rm vol}_{g_i}))}(W_{h_i}TW_{h_i}^{\ast})=\mu_{B(L_2(X,{\rm vol}_G))}(T).
\end{equation}

It is immediate that
$$W_{h_i}TW_{h_i}^{\ast}=W_{h_i}SW_{h_i}^{\ast}\cdot W_{h_i}M_{\phi}(1+\Delta_{G})^{-1}M_{\phi}W_{h_i}^{\ast}.$$
Using Lemma \ref{yuri lemma}, we write
$$W_{h_i}TW_{h_i}^{\ast}-W_{h_i}SW_{h_i}^{\ast}\cdot M_{\phi\circ h_i^{-1}}(1+\Delta_{g_i})^{-1}M_{\phi\circ h_i^{-1}}\in\mathcal{L}_{\frac{d}{3},\infty}(L_2(\mathbb{R}^d,{\rm vol}_{g_i})).$$
Let $P_{i}=M_{{\rm det}^{-\frac14}(g_{i})}\Delta_{g_i}M_{{\rm det}^{\frac14}(g_{i})}$. By Theorem \ref{euclidean compactness of commutators},
\begin{align*}
M_{{\rm det}^{-\frac14}(g_{i})}M_{\phi\circ h_i^{-1}}[(1+\Delta_{g_i})^{-1},M_{\phi\circ h_i^{-1}}]
M_{{\rm det}^{\frac14}(g_{i})}
&=M_{\phi\circ h_i^{-1}}[(1+P_{i})^{-1},M_{\phi\circ h_i^{-1}}]\\
&\in\mathcal{L}_{\frac{d}{3},\infty}(L_2(\mathbb{R}^d)),
\end{align*}
and therefore
$$M_{\phi\circ h_i^{-1}}(1+\Delta_{g_i})^{-1}M_{\phi\circ h_i^{-1}}-M_{\phi^2\circ h_i^{-1}}(1+\Delta_{g_i})^{-1}\in\mathcal{L}_{\frac{d}{3},\infty}(L_2(\mathbb{R}^d,{\rm vol}_{g_i})).$$
Since
$$W_{h_i}SW_{h_i}^{\ast}=W_{h_i}SW_{h_i}^{\ast}\cdot M_{\phi\circ h_i^{-1}},$$
it follows that
\begin{equation}\label{alp eq1}
W_{h_i}TW_{h_i}^{\ast}-W_{h_i}SW_{h_i}^{\ast}\cdot (1+\Delta_{g_i})^{-1}\in\mathcal{L}_{\frac{d}{3},\infty}(L_2(\mathbb{R}^d,{\rm vol}_{g_i})).
\end{equation}

By construction of $\Pi_X,$ we have ${\rm Ext}_{h_i(\U_i)}(W_{h_i}SW_{h_i}^{\ast})\in\Pi_{\mathbb{R}^d}.$ By Theorem \ref{for curved plane thm}, there exists a limit
\begin{align*}
&\lim_{t\to\infty}t^{\frac{2}{d}}\mu_{B(L_2(\mathbb{R}^d,{\rm vol}_{g_i}))}(t,W_{h_i}SW_{h_i}^{\ast}\cdot (1+\Delta_{g_i})^{-1})\\
&=c_d\|{\rm sym}({\rm Ext}_{h_i(\U_i)}(W_{h_i}SW_{h_i}^{\ast}))\|_{L_{\frac{d}{2}}(h_i(\U_i)\times\mathbb{R}^d,e^{-q_{g_i}}dm)}\\
&=c_d\|{\rm sym}_X(S)\|_{L_{\frac{d}{2}}(T^{\ast}X,e^{-q_G}d\lambda)}.
\end{align*}
Using Lemma \ref{bs separable lemma} and taking \eqref{alp eq1} into account, we conclude that
$$\lim_{t\to\infty}t^{\frac2d}\mu_{B(L_2(\mathbb{R}^d,{\rm vol}_{g_i}))}(t,W_{h_i}TW_{h_i}^{\ast})=c_{p,d}\|{\rm sym}_X(S)\|_{L_p(T^{\ast}X,e^{-q_G}d\lambda)}.$$
Using \eqref{alp eq0}, we write
\begin{equation}\label{alp eq2}
\lim_{t\rightarrow\infty}t^{\frac2d}\mu(t,T)
=c_d\Big\|\sym_{X}(S)\Big\|_{L_{\frac{d}{2}}(T^{\ast}X,e^{-q_G}d\lambda)}.
\end{equation}

Since $S=SM_{\phi}$ and since (from Lemma \ref{manifold compactness of commutators})
$$M_{\phi}(1+\Delta_G)^{-1}M_{\phi}-M_{\phi^2}(1+\Delta_G)^{-1}\in\mathcal{L}_{\frac{d}{3},\infty}(L_2(X,{\rm vol}_G)),$$
it follows that
$$T-S(1+\Delta_G)^{-1}\in\mathcal{L}_{\frac{d}{3},\infty}(L_2(X,{\rm vol}_G)).$$
The assertion follows now from \eqref{alp eq2} and Lemma \ref{bs separable lemma}.
\end{proof}

\section{Asymptotics on manifold}\label{sec-7}
In this section, we complete the proof of Theorem \ref{asymptotic on manifold}. The case where $p=\frac{d}{2}$ is treated in Subsection \ref{asymptotics for d2}, while the general case is addressed in Subsection \ref{proof of genral Theorem 1.2}. In this section, let $c_{d}$ given in Theorem \ref{asymptotic on Rd} in the sequel.
\subsection{Asymptotics for $p=\frac{d}{2}$}\label{asymptotics for d2}
We prove Theorem \ref{asymptotic on manifold} with $p=\frac{d}{2}$ in this subsection. To do that, we demand the following approximation on manifold.
\begin{lemma}\label{approximate lemma}
Let $X$ be a compact Riemannian manifold with a good atlas $\{(\U_{i},h_{i})\}_{i\in\mathbb{I}}$. There is a sequence $\{\phi_{n,k}\}_{n\in\mathbb{N},1\leq k\leq K}$ such that 
\begin{enumerate}[\rm(a)]
\item for all $n\in\mathbb{N}$ and $1\leq k\leq K$, function $\phi_{n,k}$ is a smooth function on $X$;
\item for all $n\in\mathbb{N}$ and $1\leq k\leq K$, function $0\leq\phi_{n,k}\leq1$;
\item\label{convergence to 1} for all $n\in\mathbb{N},1\leq k\leq K$, function $\phi_{n,k}$ is supported in some chart $(\U_{i_{k}},h_{i_{k}})$;
\item the sum $\sum_{k=1}^{K}\phi_{n,k}^{3}\rightarrow1$ in the sense of almost everywhere as $n\rightarrow\infty$;
\item\label{orthogonality of appro} for all $n\in\mathbb{N}$, the functions $\{\phi_{n,k}\}_{k=1}^{K}$ are disjoint supported.
\end{enumerate}
\end{lemma}
\begin{proof} Choose pairwise disjoint open sets $\{\Delta_k\}_{k=1}^K$ such that $\Delta_{k}$ is contained in some $\mathcal{U}_{i_{k}}$ and their union has full measure. For $1\leq k\leq K$, choose a sequence $\{\phi_{n,k}\}_{n\in\mathbb{N}}$ such that (i) $\phi_{n,k}\in C^{\infty}_c(\Delta_k)$ is positive taking values in $[0,1]$ for every $n\in\mathbb{N}$ and (ii) $\phi_{n,k}\uparrow \chi_{\Delta_k}$ almost everywhere as $n\to\infty.$ The assertions follow.
\end{proof}

\begin{lemma}\label{orthogonal op limit}
Let $\phi\in C^{\infty}(X)$ be supported in some chart $(\U_{i},h_{i})$ (which belongs to a good atlas). For any $S\in\Pi_{X}$, we have 
$$\lim_{t\rightarrow\infty}t^{\frac2d}\mu(t,M_{\phi}SM_{\phi}(1+\Delta_{G})^{-1}M_{\phi})
=c_d\|\sym_{X}(M_{\phi^{3}}S)\|_{L_{\frac{d}{2}}(T^{\ast}X,e^{-q_{G}}d\lambda)}.$$
\end{lemma}
\begin{proof} Applying Theorem \ref{asymptotic local} to the operator $M_{\phi}SM_{\phi^2},$ we conclude that
$$\lim_{t\rightarrow\infty}t^{\frac{1}{p}}\mu(t,M_{\phi}SM_{\phi^2}(1+\Delta_{G})^{-\frac{d}{2p}})
=c_{p,d}\Big\|\sym_{X}(M_{\phi}SM_{\phi^2})\Big\|_{L_p(T^{\ast}X,e^{-q_{G}}d\lambda)}.$$
	
Write 
$$M_{\phi}SM_{\phi}(1+\Delta_{G})^{-1}M_{\phi}=M_{\phi}SM_{\phi^2}(1+\Delta_{G})^{-1}
+M_{\phi}SM_{\phi}\cdot [(1+\Delta_{G})^{-1},M_{\phi}].$$
By Lemma \ref{manifold compactness of commutators}, we have
$$M_{\phi}SM_{\phi}[(1+\Delta_{G})^{-1},M_{\phi}]\in\mathcal{L}_{\frac{d}{3},\infty}(L_2(X,{\rm vol}_G)).$$
Therefore, by Lemma \ref{bs separable lemma},
$$\lim_{t\rightarrow\infty}t^{\frac2d}\mu(t,M_{\phi}SM_{\phi}(1+\Delta_{G})^{-1}M_{\phi})
=c_d\|\sym_X(M_{\phi}SM_{\phi^2})\|_{L_{\frac{d}{2}}(T^{\ast}X,e^{-q_{G}}d\lambda)}.$$
Taking into account that
\begin{align*}
\sym_X(M_{\phi}SM_{\phi^2})
&=\sym_X(M_{\phi})\cdot\sym_X(S)\cdot\sym_S(M_{\phi})^2\\
&=\sym_X(M_{\phi^3})\cdot\sym_X(S)\\
&=\sym_X(M_{\phi^3}S),
\end{align*}
we complete the proof.
\end{proof}

\begin{lemma}\label{limit approximate factor}
For each $n\in\mathbb{N}$, let $\psi_{n}=\sum_{k=1}^{K}\phi_{n,k}^{3}$ where sequence $\{\phi_{n,k}\}_{n\in\mathbb{N},1\leq k\leq K}$ is given in Lemma \ref{approximate lemma}. For any $S\in\Pi_X$, we have 
$$\lim_{t\rightarrow\infty}t^{\frac2d}\mu(t,SM_{\psi_n}(1+\Delta_{G})^{-1})
=c_d\|\sym_{X}(M_{\psi_n}S)\|_{L_{\frac{d}{2}}(T^{\ast}X,e^{-q_{G}}d\lambda)}.$$ 
\end{lemma}
\begin{proof} By Lemma \ref{orthogonal op limit}, for every $n\in\mathbb{N}$ and for every $1\leq k\leq K,$ there exists a limit
$$\lim_{t\rightarrow\infty}t^{\frac2d}\mu(t,M_{\phi_{n,k}}SM_{\phi_{n,k}}(1+\Delta_{G})^{-1}M_{\phi_{n,k}})
=c_d\|\sym_{X}(M_{\phi_{n,k}^{3}}S)\|_{L_{\frac{d}{2}}(T^{\ast}X,e^{-q_{G}}d\lambda)}.$$	
For every $n\in\mathbb{N},$ the operators $$\{M_{\phi_{n,k}}SM_{\phi_{n,k}}(1+\Delta_{G})^{-1}M_{\phi_{n,k}}\}_{1\leq k\leq K}$$
are pairwise orthogonal. By Lemma \ref{bs direct sum lemma}, there exists a limit
\begin{align*}
&\lim_{t\rightarrow\infty}t^{\frac2d}\mu(t,\sum_{k=1}^KM_{\phi_{n,k}}SM_{\phi_{n,k}}
(1+\Delta_{G})^{-1}M_{\phi_{n,k}})\\
&=c_d\Big(\sum_{k=1}^K
\|\sym_{X}(M_{\phi_{n,k}^{3}}S)\|_{L_{\frac{d}{2}}(T^{\ast}X,e^{-q_{G}}d\lambda)}^{\frac{d}{2}}\Big)^{\frac2d}.
\end{align*}
	
Write 
\begin{align*}
SM_{\psi_{n}}(1+\Delta_{G})^{-1}
&=\sum_{k=1}^{K}SM_{\phi_{n,k}^{3}}(1+\Delta_G)^{-1}\\
&=\sum_{k=1}^{K}[S,M_{\phi_{n,k}}]\cdot M_{\phi_{n,k}^{2}}(1+\Delta_{G})^{-1}+\\
&\quad+\sum_{k=1}^{K}M_{\phi_{n,k}}SM_{\phi_{n,k}}\cdot [M_{\phi_{n,k}},(1+\Delta_{G})^{-1}]+\\
&\quad+\sum_{k=1}^{K}M_{\phi_{n,k}}SM_{\phi_{n,k}}(1+\Delta_{G})^{-1}M_{\phi_{n,k}}.
\end{align*}
Note that $[S,M_{\phi_{n,k}}]$ is compact (this is part of the definition of $\Pi_X$). Hence,
$$[S,M_{\phi_{n,k}}]\cdot M_{\phi_{n,k}^{2}}(1+\Delta_{G})^{-1}\in(\mathcal{L}_{\frac{d}{2},\infty})_0(L_2(X,{\rm vol}_G)).$$
By Lemma \ref{manifold compactness of commutators}
$$M_{\phi_{n,k}}SM_{\phi_{n,k}}\cdot [M_{\phi_{n,k}},(1+\Delta_{G})^{-1}]\in(\mathcal{L}_{\frac{d}{2},\infty})_0(L_2(X,{\rm vol}_G)).$$
Hence,
$$SM_{\psi_n}(1+\Delta_{G})^{-1}-\sum_{k=1}^{K}M_{\phi_{n,k}}SM_{\phi_{n,k}}(1+\Delta_{G})^{-1}M_{\phi_{n,k}}\in (\mathcal{L}_{\frac{d}{2},\infty})_0(L_2(X,{\rm vol}_G)).$$
By Lemma \ref{bs separable lemma},
\begin{align*}
&\lim_{t\rightarrow\infty}t^{\frac2d}\mu(t,SM_{\psi_n}(1+\Delta_{G})^{-1})\\
&=c_d(\sum_{k=1}^K
\|\sym_{X}(M_{\phi_{n,k}^{3}}S)\|_{L_{\frac{d}{2}}(T^{\ast}X,e^{-q_{G}}d\lambda)}^{\frac{d}{2}})^{\frac2d}\\
&=c_d(\sum_{k=1}^K\|\sym_{X}(M_{\phi_{n,k}})^3
\cdot\sym_X(S)\|_{L_{\frac{d}{2}}(T^{\ast}X,e^{-q_{G}}d\lambda)}^{\frac{d}{2}})^{\frac2d}.
\end{align*}
Since the functions $\{\phi_{n,k}\}_{1\leq k\leq K}$ are pairwise disjointly supported, it follows that so are the functions $\{{\rm sym}_X(M_{\phi_{n,k}})\}_{1\leq k\leq K}.$ Hence,
\begin{align*}
&(\sum_{k=1}^K\|\sym_{X}(M_{\phi_{n,k}})^3\cdot\sym_X(S)\|_{L_{\frac{d}{2}}(T^{\ast}X,e^{-q_{G}}d\lambda)}^{\frac{d}{2}})^{\frac2d}\\
&=\|\sym_{X}(M_{\psi_n})\cdot\sym_X(S)\|_{L_{\frac{d}{2}}(T^{\ast}X,e^{-q_{G}}d\lambda)}.
\end{align*}
This completes the proof.
\end{proof}

We firstly provide the proof for Theorem \ref{asymptotic on manifold} in the case $p=\frac{d}{2}$.
\begin{proof}[Proof of Theorem \ref{asymptotic on manifold} for $p=\frac{d}{2}$]
For each $n\in\mathbb{N}$, let $\psi_{n}=\sum_{k=1}^{K}\phi_{n,k}^{3}$ where the sequence $\{\phi_{n,k}\}_{n\in\mathbb{N},1\leq k\leq K}$ is given in Lemma \ref{approximate lemma}. It is immediate that
\begin{align*}
&\|SM_{\psi_n}(1+\Delta_{G})^{-1}-S(1+\Delta_{G})^{-1}\|_{\mathcal{L}_{\frac{d}{2},\infty}(L_2(X,{\rm vol}_G))}\\
&\leq\|S\|_{\infty}\|M_{\psi_n-1}(1+\Delta_{G})^{-1}\|_{\mathcal{L}_{\frac{d}{2},\infty}(L_2(X,{\rm vol}_G))}.
\end{align*}
Applying Lemma \ref{special cwikel manifold}, we conclude
$$\|M_{\psi_n-1}(1+\Delta_{G})^{-1}\|_{\mathcal{L}_{\frac{d}{2},\infty}}\leq c_{X,g}\|\psi_n-1\|_{\max\{\frac{d}{2},3\}}\to0,\quad n\to\infty.$$
Here, the convergence is guaranteed by Lemma \ref{approximate lemma} and by the Dominated Convergence Theorem.
 
Lemma \ref{limit approximate factor} asserts that
$$\lim_{t\rightarrow\infty}t^{\frac2d}\mu(t,SM_{\psi_n}(1+\Delta_{G})^{-1})
=c_d\|\sym_{X}(M_{\psi_n}S)\|_{L_{\frac{d}{2}}(T^{\ast}X,e^{-q_{G}}d\lambda)}.$$ 
Using Lemma \ref{bs limit lemma}, we conclude that the following limits exist and are equal:
\begin{align*}
\lim_{t\rightarrow\infty}t^{\frac2d}\mu(t,S(1+\Delta_{G})^{-1})
&=c_d\lim_{n\to\infty}\|\sym_X(M_{\psi_n}S)\|_{L_{\frac{d}{2}}(T^{\ast}X,e^{-q_{G}}d\lambda)}\\
&=c_d\|\sym_X(S)\|_{L_{\frac{d}{2}}(T^{\ast}X,e^{-q_{G}}d\lambda)}.
\end{align*}
Here, the very last equality is guaranteed by Lemma \ref{approximate lemma} and by the Dominated Convergence Theorem. This completes the proof.
\end{proof}

\subsection{General case for an arbitrary $p$}\label{proof of genral Theorem 1.2}
In this subsection, we complete the proof of Theorem \ref{asymptotic on manifold}.
\begin{lemma}\label{Mpsi P commutator lemma}
Let $P:W^{2,2}(\mathbb{R}^d)\to L_2(\mathbb{R}^d)$ be a positive, self-adjoint, elliptic differential operator whose coefficients are smooth and also are constant outside some ball. For every $\psi\in C_0(\mathbb{R}^d)$ and for every $\phi\in C^{\infty}_c(\mathbb{R}^d)$ we have
$$[M_{\psi},M_{\phi}(1+P)^{-1}M_{\phi}]\in(\mathcal{L}_{\frac{d}{2},\infty})_0(L_2(\mathbb{R}^d)).$$ 
\end{lemma}
\begin{proof}
Write
\begin{align*}
[M_{\psi},M_{\phi}(1+P)^{-1}M_{\phi}]=[M_{\psi\phi},(1+P)^{-1}]M_{\phi}
+[(1+P)^{-1},M_{\phi}]M_{\psi\phi}.
\end{align*}
By theorem \ref{euclidean compactness of commutators}, the factors on the right hand side belong to  $\mathcal{L}_{\frac{d}{3},\infty}(L_2(\mathbb{R}^d))$ and therefore $(\mathcal{L}_{\frac{d}{2},\infty})_0(L_2(\mathbb{R}^d))$.
\end{proof}

\begin{lemma}\label{Rj P commutator lemma}
Let $P:W^{2,2}(\mathbb{R}^d)\to L_2(\mathbb{R}^d)$ be a positive, self-adjoint, elliptic differential operator whose coefficients are smooth and also are constant outside some ball. For every $\phi\in C^{\infty}_c(\mathbb{R}^d)$ we have
$$[R_j,M_{\phi}(1+P)^{-1}M_{\phi}]\in(\mathcal{L}_{\frac{d}{2},\infty})_0(L_2(\mathbb{R}^d)),\quad 1\leq j\leq d.$$ 
\end{lemma}
\begin{proof}
Write
\begin{align*}
[R_j,M_{\phi}(1+P)^{-1}M_{\phi}]=R_j M_{\phi}[(1+P)^{-1},M_{\phi}]+[(1+P)^{-1},M_{\phi}]M_{\phi}R_j+\\
+[R_j,M_{\phi}]M_{\phi}(1+P)^{-1}+(1+P)^{-1}M_{\phi}[R_j,M_{\phi}]+\\
+[M_{\phi}R_j M_{\phi},(1+P)^{-1}].
\end{align*}
By Theorem \ref{euclidean compactness of commutators} the first and second factors on the right hand side belong to  $\mathcal{L}_{\frac{d}{3},\infty}(L_2(\mathbb{R}^d))$ and therefore $(\mathcal{L}_{\frac{d}{2},\infty})_0(L_2(\mathbb{R}^d))$. 
Since operator $[R_j,M_{\phi}]$ is compact (see e.g. \cite{Uc1978,JW1982}) and $M_{\phi}(1+P)^{-1},(1+P)^{-1}M_{\phi}$ belong to $\mathcal{L}_{\frac{d}{2},\infty}(L_2(\mathbb{R}^d))$ by Lemma \ref{P-Mf compactness}, the third and fourth factors on the right hand side belong to $(\mathcal{L}_{\frac{d}{2},\infty})_0(L_2(\mathbb{R}^d))$. 
For the fifth factor on the right hand side, write 
\begin{align*}
&[M_{\phi}R_{j}M_{\phi},(1+P)^{-1}]=\\
&=(1+P)^{-\frac{1}{2}}[M_{\phi}R_{j}M_{\phi},(1+P)^{-\frac{1}{2}}]
+[M_{\phi}R_{j}M_{\phi},(1+P)^{-\frac{1}{2}}](1+P)^{-\frac{1}{2}}.
\end{align*}
Applying Theorem \ref{1 over 2 case}, it gives that operator $[M_{\phi}R_{j}M_{\phi},(1+P)^{-1}]$ belongs to $\mathcal{L}_{\frac{d}{3},\infty}(L_2(\mathbb{R}^d))$ and therefore $(\mathcal{L}_{\frac{d}{2},\infty})_{0}(L_2(\mathbb{R}^d))$. The proof is complete.
\end{proof}

\begin{lemma}\label{S P commutator lemma} Let $P:W^{2,2}(\mathbb{R}^d)\to L_2(\mathbb{R}^d)$ be a positive, self-adjoint, elliptic differential operator whose coefficients are smooth and also are constant outside some ball. For every $S\in\Pi_{\mathbb{R}^d}$ and for every $\phi\in C^{\infty}_c(\mathbb{R}^d)$ we have
$$[S,M_{\phi}(1+P)^{-1}M_{\phi}]\in(\mathcal{L}_{\frac{d}{2},\infty})_0(L_2(\mathbb{R}^d)).$$ 
\end{lemma}
\begin{proof} Let $S$ belong to $\ast$-algebra generated by $\{M_f:\ f\in C_0(\mathbb{R}^d)\}$ and by $\{R_j\}_{j=1}^d.$ Using Lemmas \ref{Mpsi P commutator lemma}, \ref{Rj P commutator lemma} and the Leibniz rule, we conclude that
$$[S,M_{\phi}(1+P)^{-1}M_{\phi}]\in(\mathcal{L}_{\frac{d}{2},\infty})_0(L_2(\mathbb{R}^d)).$$ 

Consider now the general case. By the definition of $\Pi_{\mathbb{R}^d},$ we can find a sequence $\{S_n\}_{n\geq0}$ such that each $S_n$ belongs to $\ast$-algebra generated by $\{M_f:\ f\in C_0(\mathbb{R}^d)\}$ and $\{R_j\}_{j=1}^d$, and such that $S_n\to S$ in the uniform norm. By the preceding paragraph,
$$[S_n,M_{\phi}(1+P)^{-1}M_{\phi}]\in(\mathcal{L}_{\frac{d}{2},\infty})_0(L_2(\mathbb{R}^d)),\quad n\geq0.$$
We clearly have
$$[S_n,M_{\phi}(1+P)^{-1}M_{\phi}]\to [S,M_{\phi}(1+P)^{-1}M_{\phi}],\quad n\to\infty.$$ 
Since $(\mathcal{L}_{\frac{d}{2},\infty})_0(L_2(\mathbb{R}^d))$ is closed in $\mathcal{L}_{\frac{d}{2},\infty}(L_2(\mathbb{R}^d)),$ the assertion follows.
\end{proof}

\begin{lemma}\label{S curved plane commutator lemma} Let $g:\mathbb{R}^d\to{\rm GL}_+(d,\mathbb{R})$ be smooth mapping constant outside some ball. For every $S\in\Pi_{\mathbb{R}^d}$ and for every $\phi\in C^{\infty}_c(\mathbb{R}^d),$ we have
$$[S,M_{\phi}(1+\Delta_g)^{-1}M_{\phi}]\in(\mathcal{L}_{\frac{d}{2},\infty})_0(L_2(\mathbb{R}^d,{\rm vol}_g)).$$ 
\end{lemma}
\begin{proof} Consider the unitary mapping $M_{{\rm det}^{\frac14}(g)}:L_2(\mathbb{R}^d,{\rm vol}_g)\to L_2(\mathbb{R}^d).$ We have
\begin{align*}
M_{{\rm det}^{-\frac14}(g)}[S,M_{\phi}(1+\Delta_g)^{-1}M_{\phi}]M_{{\rm det}^{\frac14}(g)}
=[M_{{\rm det}^{-\frac14}(g)}SM_{{\rm det}^{\frac14}(g)},M_{\phi}(1+P)^{-1}M_{\phi}]
\end{align*}
where
$$P=M_{{\rm det}^{-\frac14}(g)}\Delta_g M_{{\rm det}^{\frac14}(g)}.$$
Since $P:W^{2,2}(\mathbb{R}^d)\to L_2(\mathbb{R}^d)$ be a positive, self-adjoint, elliptic differential operator whose coefficients are smooth and also are constant outside some ball, it follows from Lemma \ref{S P commutator lemma} that
$$[M_{{\rm det}^{-\frac14}(g)}SM_{{\rm det}^{\frac14}(g)},M_{\phi}(1+P)^{-1}M_{\phi}]\in(\mathcal{L}_{\frac{d}{2},\infty})_0(L_2(\mathbb{R}^d)).$$
This suffices to complete the proof.
\end{proof}

\begin{lemma}\label{pix localised commutator lemma} Let $(X,G)$ be a compact Riemannian manifold and let $\Delta_G$ be the Laplace-Beltrami operator. Let $(\U_i,h_i)$ be a chart in a good atlas. For every $S\in\Pi_X$ compactly supported in $\U_i,$ we have
$$[S,(1+\Delta_G)^{-1}]\in(\mathcal{L}_{\frac{d}{2},\infty})_0(L_2(X,{\rm vol}_G)).$$
\end{lemma}
\begin{proof} Fix $\phi\in C^{\infty}_c(\U_i)$ such that $S=M_{\phi}S=SM_{\phi}.$ We have
\begin{align*}
[S,(1+\Delta_G)^{-1}]
&=S(1+\Delta_G)^{-1}-(1+\Delta_G)^{-1}S\\
&=SM_{\phi^2}(1+\Delta_G)^{-1}-(1+\Delta_G)^{-1}M_{\phi^2}S\\
&=[S,M_{\phi}(1+\Delta_G)^{-1}M_{\phi}]
+S\cdot[M_{\phi},(1+\Delta_G)^{-1}]+\\
&\quad+[M_{\phi},(1+\Delta_G)^{-1}]\cdot S.
\end{align*}
By Lemma \ref{manifold compactness of commutators},
$$[M_{\phi},(1+\Delta_G)^{-1}]\in\mathcal{L}_{\frac{d}{3},\infty}(L_2(X,{\rm vol}_G)).$$
Thus,
$$[S,(1+\Delta_G)^{-1}]-[S,M_{\phi}(1+\Delta_G)^{-1}M_{\phi}]\in\mathcal{L}_{\frac{d}{3},\infty}(L_2(X,{\rm vol}_G)).$$
Lemma \ref{yuri lemma} asserts that
$$M_{\phi}(1+\Delta_G)^{-1}M_{\phi}-W_{h_i}^{\ast}M_{\phi\circ h_i^{-1}}(1+\Delta_{g_i})^{-1}M_{\phi\circ h_i^{-1}}W_{h_i}\in\mathcal{L}_{\frac{d}{3},\infty}(L_2(X,{\rm vol}_G)).$$
Thus,
$$[S,(1+\Delta_G)^{-1}]-[S,W_{h_i}^{\ast}M_{\phi\circ h_i^{-1}}(1+\Delta_{g_i})^{-1}M_{\phi\circ h_i^{-1}}W_{h_i}]\in\mathcal{L}_{\frac{d}{3},\infty}(L_2(X,{\rm vol}_G)).$$
In other words,
$$[S,(1+\Delta_G)^{-1}]-W_{h_i}^{\ast}[S_i,M_{\phi\circ h_i^{-1}}(1+\Delta_{g_i})^{-1}M_{\phi\circ h_i^{-1}}]W_{h_i}\in\mathcal{L}_{\frac{d}{3},\infty}(L_2(X,{\rm vol}_G)),$$
where $S_i=W_{h_i}SW_{h_i}^{\ast}\in\Pi_{\mathbb{R}^d}.$ The assertion follows now from Lemma \ref{S curved plane commutator lemma}.
\end{proof}

\begin{lemma}\label{pix commutator lemma} Let $(X,G)$ be a compact Riemannian manifold and let $\Delta_G$ be the Laplace-Beltrami operator. For every $S\in\Pi_X,$ we have
$$[S,(1+\Delta_G)^{-1}]\in(\mathcal{L}_{\frac{d}{2},\infty})_0(L_2(X,{\rm vol}_G)).$$
\end{lemma}
\begin{proof} Let $\{(U_n,h_n)\}_{n=1}^N$ be a good atlas and let $\{\phi_k\}_{k=1}^K$ be a partition of unity such that, for every $1\leq k_1,k_2\leq K,$ either there exists $n_{k_1,k_2}$ such that $\phi_{k_1},\phi_{k_2}\in C_c(\U_{n_{k_1,k_2}})$ or $\phi_{k_1}\phi_{k_2}\equiv0.$ We write
$$S=\sum_{k_1,k_2=1}^KM_{\phi_{k_1}}SM_{\phi_{k_2}}.$$
If there exists $n_{k_1,k_2}$ such that $\phi_{k_1},\phi_{k_2}\in C^{\infty}_c(\U_{n_{k_1,k_2}}),$ then 
$$[M_{\phi_{k_1}}SM_{\phi_{k_2}},(1+\Delta_G)^{-1}]\in(\mathcal{L}_{\frac{d}{2},\infty})_0(L_2(X,{\rm vol}_G))$$
by Lemma \ref{pix localised commutator lemma}. If $\phi_{k_1}\phi_{k_2}=0,$ then
$$M_{\phi_{k_1}}SM_{\phi_{k_2}}=M_{\phi_{k_1}}\cdot [S,M_{\phi_{k_2}}]\in\mathcal{K}(L_2(X,{\rm vol}_G)).$$
In either case, we have
$$[M_{\phi_{k_1}}SM_{\phi_{k_2}},(1+\Delta_G)^{-1}]\in(\mathcal{L}_{\frac{d}{2},\infty})_0(L_2(X,{\rm vol}_G)).$$
Summing over $1\leq k_1,k_2\leq K,$ we complete the proof.
\end{proof}

Now, let us give a complete proof for Theorem \ref{asymptotic on manifold}.
\begin{proof}[Proof of Theorem \ref{asymptotic on manifold}] Without loss of generality, $S\geq0.$ Using Theorem \ref{asymptotic on manifold} for the case $p=\frac{d}{2},$ we write
$$\lim_{t\rightarrow\infty}t^{\frac2d}\mu(t,S^{\frac{2p}{d}}(1+\Delta_{G})^{-1})=c_d\|\sym_X(S^{\frac{2p}{d}})\|_{L_{\frac{d}{2}}(T^{\ast}X,e^{-q_{G}}d\lambda)}.$$
According to Lemma \ref{pix commutator lemma},
$$S^{\frac{2p}{d}}(1+\Delta_{G})^{-1}-S^{\frac{p}{d}}(1+\Delta_G)^{-1}S^{\frac{p}{d}}\in(\mathcal{L}_{\frac{d}{2},\infty})_0(L_2(X,{\rm vol}_G)).$$
It follows from Lemma \ref{bs separable lemma} that
$$\lim_{t\rightarrow\infty}t^{\frac2d}\mu(t,S^{\frac{p}{d}}(1+\Delta_G)^{-1}S^{\frac{p}{d}})=c_d\|\sym_X(S^{\frac{2p}{d}})\|_{L_{\frac{d}{2}}(T^{\ast}X,e^{-q_{G}}d\lambda)}.$$
Taking this equality to the power $\frac{d}{2p},$ we write
$$\lim_{t\rightarrow\infty}t^{\frac1p}\mu\Big(t,\Big(S^{\frac{p}{d}}(1+\Delta_G)^{-1}S^{\frac{p}{d}}\Big)^{\frac{d}{2p}}\Big)=c_d^{\frac{d}{2p}}\|\sym_X(S^{\frac{2p}{d}})\|_{L_{\frac{d}{2}}(T^{\ast}X,e^{-q_{G}}d\lambda)}^{\frac{d}{2p}}.$$
Applying Theorem \ref{main example thm} with $s=\frac{d}{2p}$ and with $S^{\frac{2p}{d}}$ instead of $S,$ we write
$$S(1+\Delta_G)^{-\frac{d}{2p}}-(S^{\frac{p}{d}}(1+\Delta_G)^{-1}S^{\frac{p}{d}})^{\frac{d}{2p}}\in(\mathcal{L}_{p,\infty})_0(L_2(X,{\rm vol}_G)).$$
It follows from Lemma \ref{bs separable lemma} that
$$\lim_{t\rightarrow\infty}t^{\frac1p}\mu\Big(t,S(1+\Delta_G)^{-\frac{d}{2p}}\Big)=c_d^{\frac{d}{2p}}\|\sym_X(S^{\frac{2p}{d}})\|_{L_{\frac{d}{2}}(T^{\ast}X,e^{-q_{G}}d\lambda)}^{\frac{d}{2p}}.$$
Since ${\rm sym}_X$ is an $\ast$-homomorphism, it follows that
$$\sym_X(S^{\frac{2p}{d}})=\sym_X(S)^{\frac{2p}{d}}$$
and the assertion follows. The proof is complete.
\end{proof}

{\bf Acknowledgement:} The first and the third author are supported by the ARC DP230100434. The second author is supported by NSFC No. 12371138.

\end{document}